\documentclass[twoside,10pt]{article}

\usepackage[pdftex]{graphicx,xcolor}     
\usepackage[]{hyperref}                        
\usepackage{amsmath}
\usepackage{amsthm}
\usepackage{bm,latexsym}
\usepackage{amsfonts}
\usepackage{indentfirst}
\usepackage{amsxtra}
\usepackage{amssymb}
\newtheorem{Lemma}{Lemma}[section]
\newtheorem{Theorem}{Theorem}

\newtheorem{Proposition}{Proposition}[section]

\theoremstyle{remark}
\newtheorem{Remark}{Remark}[section]
\numberwithin{equation}{section}

\begin{document}

\title{\bf
Viscous Shock Wave to an Inflow Problem for Compressible Viscous Gas with Large Density Oscillations}
\author{{\bf Dongfen Bian}\\
School of Mathematics and Statistics, Beijing Institute of Technology\\
Beijing 100081, China\\
and\\
Department of Mathematics, City University of Hong Kong, Hong Kong, China\\
E-mail address: biandongfen@bit.edu.cn\\[2mm]
{\bf  Lili Fan}\\
School of Mathematics and Computer Science,
Wuhan Polytechnic University\\
 Wuhan 430023, China\\
 E-mail address: fll810@live.cn\\[2mm]
{\bf Lin He}\\
School of Mathematics and Statistics, Wuhan University\\
Wuhan 430072, China\\
E-mail address: linhe1989@whu.edu.cn\\[2mm]
{\bf Huijiang Zhao}\\
School of Mathematics and Statistics, Wuhan University\\
Wuhan 430072, China\\
and\\
Computational Science Hubei Key Laboratory, Wuhan University\\
Wuhan 430072, China\\
E-mail address: hhjjzhao@hotmail.com}

\date{} \maketitle


\begin{abstract}
This paper is concerned with the inflow problem for the one-dimensional compressible Navier-Stokes equations. For such a problem, F. M. Huang, A. Matsumura and X. D. Shi showed in \cite{Huang-Matsumura-Shi-CMP-2003} that there exists viscous shock wave solution to the inflow problem and both the boundary layer solution, the viscous shock wave, and their superposition are time-asymptotically nonlinear stable under small initial perturbation. The main purpose of this paper is to show that similar stability results still hold for a class of large initial perturbation which can allow the initial density to have large oscillations. The proofs are given by an elementary energy method and the key point is to deduce the desired uniform positive lower and upper bounds on the density.

\bigbreak
\end{abstract}

\tableofcontents

\section{Introduction}
\subsection{The inflow problem}
This paper is concerned with the large time behaviors of solutions to the inflow problem for one-dimensional compressible Navier-Stokes equations in the Eulerian coordinates
\begin{equation} \label{EulerI}
  \left\{
  \begin{array}{rl}
    \rho_\tau+(\rho u)_{\tilde{x}}&=0,\\
    (\rho u)_\tau+(\rho u^2 +p)_{\tilde{x}}&=\mu u_{\tilde{x}\tilde{x}}
  \end{array}
  \right.
\end{equation}
on the half line $\mathbb{R}_+=[0,+\infty)$ with prescribed initial and boundary conditions
\begin{equation}\label{IBC-Euler}
   \begin{cases}
    \left(\rho(\tau,\tilde{x}),u(\tau,\tilde{x})\right)|_{\tilde{x}=0}
    =(\rho_-,u_-),&u_->0,\ \rho_->0,\ \tau\geq 0,\\
    \left(\rho(\tau,\tilde{x}),u(\tau,\tilde{x})\right)|_{\tau=0}
    =\left(\rho_0\left(\tilde{x}\right),u_0\left(\tilde{x}\right)\right)\to \left(\rho_+,u_+\right), \quad &\textrm{as}\ 0\leq\tilde{x}\to+\infty,
  \end{cases}
\end{equation}
which are assumed to satisfy the compatibility condition
$$
\rho_0(0)=\rho_-,\quad u_0(0)=u_-.
$$
Here $\tilde{x}$ and $\tau$ represent the Eulerian space variable and the time variable, respectively, $\rho(\tau,\tilde{x})(>0)$, $u(\tau,\tilde{x})$, and  $p=p(\rho) =\rho^{\gamma}$ with $\gamma\geq 1$ being the adiabatic exponent are, respectively, the density, the velocity, and the pressure, while the viscosity coefficient $\mu(>0)$, $\rho_{\pm}(> 0)$, and $u_{\pm}$ are constants.

For such an initial boundary problem, as classified in \cite{Matsumura-MAA-2001}, the assumption that the boundary velocity $u_->0$ implies that the fluid with the density $\rho_-$ flow into the region $\mathbb{R}_+$ through the boundary $\tilde{x}=0$, and thus the problem \eqref{EulerI} is called the inflow problem. The cases $u_-=0$ and $u_-<0$, where the condition $\rho|_{\tilde{x}=0}=\rho_-$
is removed, are called impermeable wall problem and the outflow problem, respectively. Throughout this manuscript, we will concerned with the inflow problem \eqref{EulerI}. For the corresponding impermeable wall problem and outflow problem, those interested are referred to \cite{Matsumura-Mei-ARMA-1999, Matsumura-Nishihara-QAM-2000} and \cite{Huang-Qin-JDE-2009, Kagei-Kawashima-CMP-2006, Kawashima-Nakamura-Nishibata-Zhu-M3AS-2010, Kawashima-Nishibata-Zhu-CMP-2003, Kawahima-Zhu-JDE-2008, Nakamura-Nishibata-SIMA-2009, Nakamura-Nishibata-Yuge-JDE-2007, Zhu-report} and the references cited therein, respectively.

\subsection{The classifications of the large behaviors}
To explain the main purpose of this manuscript, we first reformulate the initial-boundary value problem \eqref{EulerI}-\eqref{IBC-Euler} as in \cite{Matsumura-Nishihara-CMP-2001}: Let $x$ be the Lagrangian space variable, $t$ be the time variable, and $v=\frac 1\rho$ denote the specific volume, we can then transform the initial-boundary value problem \eqref{EulerI}-\eqref{IBC-Euler} into the following problem in the Lagrangian coordinates:
\begin{equation} \label{LagrangeI}
  \begin{cases}
    v_t-u_x=0,   &x>s_- t,\ t>0, \\
    u_t+p(v)_x=\mu\left(\frac{u_x}{v}\right)_x, & x>s_- t,\ t>0, \\
    (v(t,x),u(t,x))|_{x=s_-t}=(v_-,u_-),&u_->0,\\
    (v(t,x),u(t,x))|_{t=0}=(v_0(x),u_0(x))\to (v_+,u_+),&\textrm{as}\ 0\leq x\to +\infty,
  \end{cases}
\end{equation}
where
\begin{equation}\label{1.4}
  p(v)=v^{-\gamma},\quad  v_{\pm}=\frac{1}{\rho_{\pm}}>0,\quad s_-=-\frac{u_-}{v_-}<0.
\end{equation}
The characteristic speeds of the corresponding hyperbolic system of~(\ref{LagrangeI}) are
\begin{equation}\label{1.5}
  \lambda_1(v)=-\sqrt{-p'(v)},\quad \lambda_2(v)=\sqrt{-p'(v)},
\end{equation}
respectively and the sound speed $c(v)$ is defined by
\begin{equation}\label{1.6}
  c(v)=v\sqrt{-p'(v)}=\sqrt{\gamma}v^{-\frac{\gamma-1}{2}}.
\end{equation}

By comparing the fluid velocity $|u|$ with the sound speed $c(v)$, one can divide the phase space $\mathbb{R}_+\times\mathbb{R}_+$ into three regions:
\begin{equation*}
  \begin{aligned}
    \Omega_{\textrm {sub}}&:=\left\{(v,u);\ |u|<c(v),v>0,u>0\right\},\\
    \Gamma_{\textrm{trans}}&:=\left\{(v,u);\ |u|=c(v),v>0,u>0\right\},\\
    \Omega_{\textrm{super}}&:=\left\{(v,u);\ |u|>c(v),v>0,u>0\right\},
  \end{aligned}
\end{equation*}
which are called the subsonic, transonic and supersonic regions, respectively.

It is now well-understood that the large time behaviors of global solutions to the Cauchy problem of the one-dimensional compressible Navier-Stokes equations \eqref{EulerI} can be described by the $i-$rarefaction wave $(V_i^{\textrm{RW}}(x/t;w_l,w_r), U_i^{\textrm{RW}}(x/t;w_l,w_r))$ ($i=1,2$) which is the unique rarefaction wave solution of the Riemann problem of the resulting Euler equations connecting the two states $w_l=(v_l,u_l)$ and $w_r=(v_r,u_r)$, the suitably shifted $i-$viscous shock wave $(V_i^{\textrm{VSW}}(x-s_it+\sigma_i;w_l,w_r), U_i^{\textrm{VSW}}(x-s_it+\sigma_i;w_l,w_r))$ ($i=1,2$) connecting $w_l$ and $w_r$
and their superpositions (For some progress on the mathematical justifications of such an expectation, see \cite{Liu-Zeng-MemoirsAMS-2015, Matsumura-Nishihara-JJAM-1985, Matsumura-Nishihara-CMP-1992, Matsumura-Nishihara-QAM-2000, Wang-Zhao-Zou-KRM-2013} and the references cited therein.) Here the $i-$viscous shock wave $(V_i^{\textrm{VSW}}(\xi;w_l,w_r), U_i^{\textrm{VSW}}(\xi;w_l,w_r))$ ($i=1,2$) is a traveling wave solution of the one-dimensional compressible Navier-Stokes equations \eqref{LagrangeI}$_1$-\eqref{LagrangeI}$_2$ connecting $w_l$ and $w_r$
which solves
\begin{equation} \label{Viscous-shock-wave}
  \left\{
  \begin{array}{rl}
   U^{\textrm{VSW}}_i(\xi;w_l,w_r)-u_l&=-s_i\left(V_i^{\textrm{VSW}}(\xi;w_l,w_r)
   -v_l\right), \\[2mm]
   \frac{s_i\mu}{V_i^{\textrm{VSW}}(\xi;w_l,w_r)}\frac{d V_i^{\textrm{VSW}}(\xi;w_l,w_r)}{d\xi}&=-s_i^2\left(V_i^{\textrm{VSW}}
   (\xi;w_l,w_r)-v_l\right)
   -\left(p\left(V_i^{\textrm{VSW}}(\xi;w_l,w_r)\right)-p(v_l)\right)\\[2mm]
   &:=h\left(V_i^{\textrm{VSW}}(\xi;w_l,w_r)\right),\\[2mm]
   V^{\textrm{VSW}}_i(-\infty;w_l,w_r)&=v_l,\quad V^{\textrm{VSW}}_i(+\infty;w_l,w_r)=v_r,
  \end{array}
  \right.
\end{equation}
where
\begin{equation}\label{wave speed}
s_i\equiv s_i(v_l,v_r)=(-1)^i\sqrt{\frac{p(v_r)-p(v_l)}{v_l-v_r}},\quad i=1,2,
\end{equation}
and the entropy condition
\begin{equation}\label{entropy condition}
u_r< u_l
\end{equation}
is assumed to be hold.

But for the initial-boundary value problem \eqref{LagrangeI}, as pointed out in \cite{Matsumura-MAA-2001}, the problem becomes complicated and to describe its large time behaviors, a new kind of nonlinear wave, the so-called {\it boundary layer solution}, or {\it BL-solution} simply in the rest of this manuscript, should be introduced which is due to the presence of the boundary. In fact, as shown in \cite{Huang-Matsumura-Shi-CMP-2003}, when $(v_-, u_-)\in \Omega_{\textrm{sub}}$, since the first wave speed $\lambda_1(v_-)$ is less than the boundary speed $s_-$, there exists a BL-solution $(V^{\textrm{BL}}(x-s_-t;w_-,w_+), U^{\textrm{BL}}(x-s_-t;w_-,w_+))$ of the one-dimensional compressible Navier-Stokes equations \eqref{LagrangeI}$_1$-\eqref{LagrangeI}$_2$ connecting $w_-=(v_-,u_-)\in \textrm{BL}(v_-,u_-)$ and some $w_+=(v_+,u_+)\in \textrm{BL}(v_-,u_-)$. Here
$$
\textrm{BL}(v_-,u_-)=\left\{(v,u)\in \Omega_{\textrm{sub}}\cup \Gamma_{\textrm{Trans}}\ \ \left|\ \ \frac uv=\frac{u_-}{v_-}=-s_-\right.\right\}
$$
denotes the BL-solution line through $(v_-,u_-)$ and $(V^{\textrm{BL}}(\xi;w_-,w_+), U^{\textrm{BL}}(\xi;w_-,w_+))$ satisfy
$$
\left\{
\begin{array}{rl}
\mu\frac{dV^{\textrm{BL}}(\xi;w_-,w_+)}{d\xi}
=&\frac{V^{\textrm{BL}}(\xi;w_-,w_+)}{s_-}
\Big\{-s^2_-\left(V^{\textrm{BL}}(\xi;w_-,w_+)-v_+\right)\\
&-\left(p\left(V^{\textrm{BL}}(\xi;w_-,w_+)\right)
-p(v_+)\right)\Big\},\\
V^{\textrm{BL}}(0;w_-,w_+)=&v_-,\quad V^{\textrm{BL}}(+\infty;w_-,w_+)=v_+,\\
U^{\textrm{BL}}(\xi;w_-,w_+)=&-s_-V^{\textrm{BL}}(\xi;w_-,w_+).
\end{array}
\right.
$$

On the other hand, since $0 >\lambda_1(v)>s_-$ in $\Omega_{\textrm{super}}$, the two characteristic fields are away from the moving boundary and hence the large time behaviors of solutions are expected to be the same as those for
the Cauchy problem. Moreover, setting
$$
\Gamma_{\textrm{Trans}}\cap \textrm{BL}(v_-,u_-)=\left\{(v_*,u_*)\right\}
$$
and noticing that $c'(v_*)>-\lambda_2(v_*)$ holds for $1<\gamma<3$,
A. Matsumura \cite{Matsumura-MAA-2001} classified all possible large time behaviors of the solutions of the initial-boundary value problem \eqref{EulerI}-\eqref{IBC-Euler} in terms of the boundary values $(v_-, u_-)$ and the far field $(v_+, u_+)$ of the initial data $(v_0(x),u_0(x))$ and it was shown in \cite{Matsumura-MAA-2001} that the large time behaviors to be expected divide the $(v, u)-$space as in Figure 1 which is taken from \cite{Matsumura-Nishihara-CMP-2001}.
\begin{figure}
\centering
\includegraphics[height=65mm]{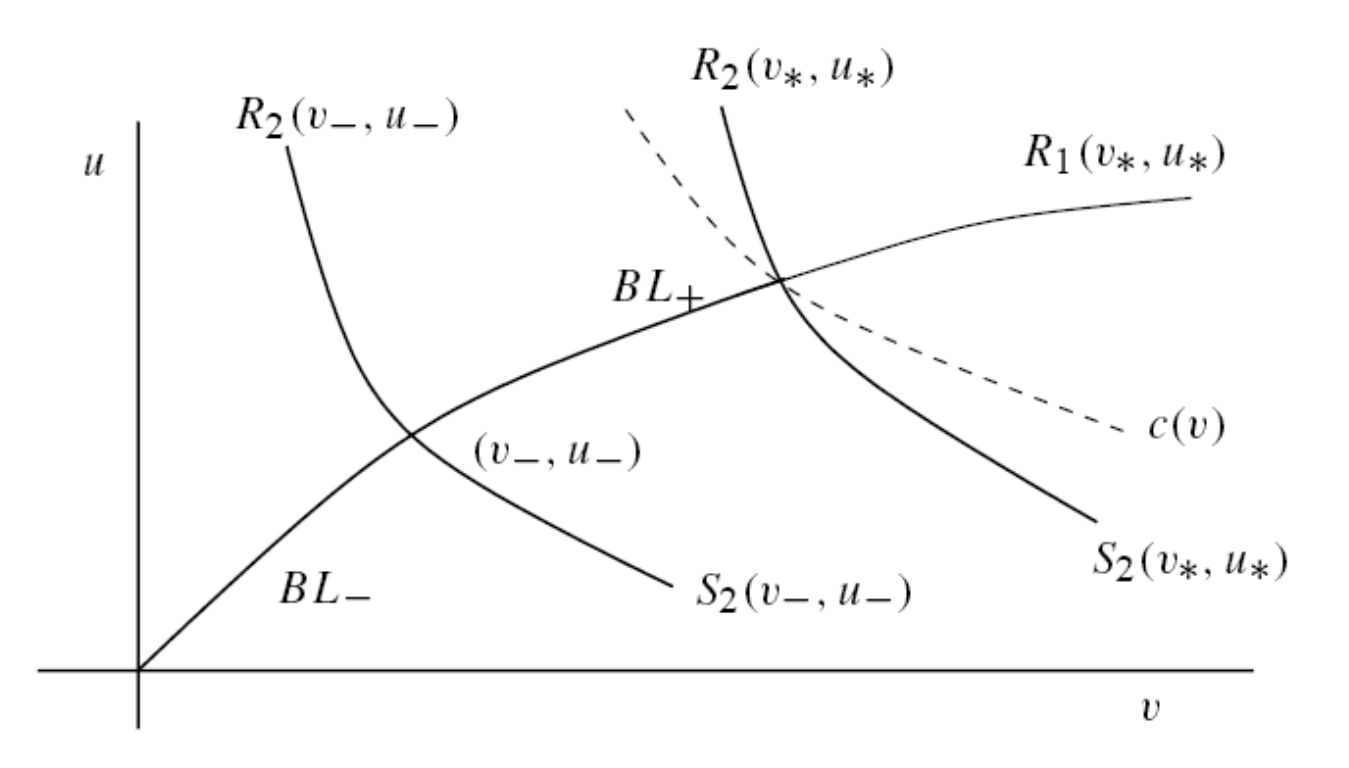}
\caption{$(v_-, u_-)\in \Omega_{\textrm{sub}}$}
\end{figure}

Here,
\begin{eqnarray*}
\textrm{BL}_+(v_-,u_-)&=&\left\{(v,u)\in \textrm{BL}(v_-,u_-):\ v_-<v\leq v_*\right\},\\
\textrm{BL}_-(v_-,u_-)&=&\left\{(v,u)\in \textrm{BL}(v_-,u_-):\ 0<v< v_-\right\},\\
\textrm{R}_1(v_*,u_*)&=&\left\{(v,u):\ u=u_*-\int^v_{v_*}\lambda_1(s)ds,\ v>v_*\right\},\\
\textrm{R}_2(v_*,u_*)&=&\left\{(v,u):\ u=u_*-\int^v_{v_*}\lambda_2(s)ds,\ v<v_*\right\},\\
\textrm{S}_2(v_-,u_-)&=&\{(v,u):\ u=u_--s_2(v_-,v)(v-v_-),\ v>v_-\},\\
\textrm{S}_2(v_*,u_*)&=&\{(v,u):\ u=u_*-s_2(v_*,v)(v-v_*),\ v>v_*\},
\end{eqnarray*}
where $s_2(a,b)$ is defined by \eqref{wave speed}.

\subsection{Former results on the mathematical justification of the expected large time behaviors and the main purpose of this manuscript}

For the mathematical justification of the above expected large time behaviors of the global solutions to the inflow problem \eqref{LagrangeI} classified in \cite{Matsumura-MAA-2001}, there are some pregoess which can be outlined as in the following:
\begin{itemize}
\item A. Matsumura and K. Nishihara \cite{Matsumura-Nishihara-CMP-2001} established the asymptotic stability of the BL-solution and the superposition of a BL-solution and a rarefaction wave for the inflow problem \eqref{LagrangeI} when $(v_-,u_-)\in \Omega_{\textrm{sub}}$,
    $(v_+,u_+)\in \textrm{BL}(v_-,u_-)\cup \textrm{BLR}_2(v_-,u_-)$ or $(v_*,u_*)\in \textrm{BL}(v_-,u_-), (v_+,u_+)\in \textrm{R}_1\textrm{R}_2(v_*,u_*)$. X. D. Shi \cite{Shi-AMASES-2003} studied the rarefaction wave case when $(v_-,u_-)\in \Omega_{\textrm{super}}, (v_+, u_+)\in \Omega_{\textrm{super}}$ and $(v_+,u_+)\in \textrm{R}_1\textrm{R}_2(v_-,u_-)$;
\item For the result concerning the viscous shock wave for the inflow problem \eqref{LagrangeI}, F. M. Huang, A. Matsumura and X. D. Shi showed in \cite{Huang-Matsumura-Shi-CMP-2003} that the viscous shock wave and the superposition of a BL-solution and a viscous shock wave for the inflow problem \eqref{LagrangeI} are nonlinear stable when $(v_-,u_-)\in \Omega_{\textrm{sub}}$, $(v_+,u_+)\in \textrm{S}_2(v_-,u_-)\cup \textrm{BLS}_2(v_-,u_-)$.
\end{itemize}

It is worth to pointing out that all the above nonlinear stability results ask that the initial perturbations are sufficiently small and strengths of some types of involving nonlinear waves such as the monotonic decreasing BL-solution, the rarefaction wave are assumed to be sufficiently small, while the strength of the monotonic increasing BL-solution is not necessarily weak and for the cases when the viscous shock wave $V_2^{\textrm{VSW}}(\xi;w_l,w_r)$ connecting the states $w_l=(v_l,u_l)$ and $w_r=(v_r,u_r)$ is involved, its strength $|v_r-v_l|$ is assumed to satisfy the following Nishida-Smoller type condition:
\begin{equation}\label{strength-VSW}
(\gamma-1)^2(v_r-v_l)<2\gamma v_l.
\end{equation}

Thus a natural question is: {\it Do similar nonlinear stability results hold for large initial perturbations?}

For Cauchy problem of the one-dimensional compressible Navier-Stokes equations \eqref{LagrangeI}$_1$-\eqref{LagrangeI}$_2$, the main difficulty lies in how to control the possible growth of its solutions induced by the nonlinearity of equations under consideration and key point is to deduce the desired uniform positive lower and upper bounds on the specific volume $v(t,x)$. For results in this direction for the case of the Cauchy problem, it is shown in \cite{Duan-Liu-Zhao-TAMS-2009, Matsumura-Nishihara-CMP-1992, Matsumura-Nishihara-QAM-2000}  that the rarefaction wave is nonlinear stable for any large initial perturbation and the nonlinear stability of viscous shock wave for a class of large initial perturbation which can allow the initial density to have large oscillations is obtained in \cite{Wang-Zhao-Zou-KRM-2013}.

For the inflow problem \eqref{LagrangeI}, in addition to the difficulty mentioned above, another difficulty is how to bound the boundary term induced by the inflow boundary condition \eqref{LagrangeI}$_3$. To overcome such a new difficulty, an argument which is based on Y. Kanel's method \cite{Kanel-Diff Eqns-1968} is introduced in \cite{Fan-Liu-Wang-Zhao-JDE-2014} for the case when its large time behavior is described by the BL-solution, the rarefaction wave, and the superposition of a BL-solution and rarefaction waves to yield
the nonlinear stability of these elementary waves for a class of large initial perturbation which can allow the initial specific volume $v_0(x)$ to have large oscillations. Moreover it is also shown in \cite{Fan-Liu-Wang-Zhao-JDE-2014} that the supersonic rarefaction wave is nonlinear stable for general large initial perturbation. Thus the results obtained in \cite{Fan-Liu-Wang-Zhao-JDE-2014} generalize the nonlinear stability results obtained in \cite{Matsumura-Nishihara-CMP-2001} and \cite{Shi-AMASES-2003} under small initial perturbations to the case of a class of large initial perturbations. Even so, no result has been obtained for the nonlinear stability of the viscous shock wave and the superposition of a BL-solution and a viscous shock wave under large initial perturbation and the main purpose of our present paper is to show that some nonlinear stability results similar to those obtained in \cite{Huang-Matsumura-Shi-CMP-2003}
still hold for a class of large initial perturbation which can allow the initial density to have large oscillations.

To simplify the presentation, we will only concentrate on the case when
$(v_-,u_-)\in \Omega_{\textrm{sub}}$, $(v_+,u_+)\in \textrm{S}_2(v_-,u_-)$ in the rest of this paper. In fact by combining the argument used in this paper with those employed in \cite{Fan-Liu-Wang-Zhao-JDE-2014}, similar result can also be obtained for the case of $(v_-,u_-)\in \Omega_{\textrm{sub}}$, $(v_+,u_+)\in \textrm{BLS}_2(v_-,u_-)$.

\subsection{Notations}
Throughout this paper, $\delta:=|v_+-v_-|$ denotes the strength of the 2$-$viscous shock wave and a positive constant $C$ is said to be $\delta-$independnet means that there exists a positive constant $C_1\geq 1$ which does not depend on $\delta$ such that $C_1^{-1}\leq C\leq C_1$, $c$, $C$ and $O(1)$ represent some $\delta-$independnet positive constant (generally large), $\epsilon$, $\lambda$ stand for some $\delta-$independent positive constants (generally small), and $C(\cdot,\cdot)$ denotes for some generic positive constant depending only on the quantities listed in the parenthesis. Notice that all the constants $c$, $C$, $C(\cdot,\cdot)$, $\epsilon$, and $\lambda$ may take different values in different places.

$A\lesssim B$ means that there is a generic $\delta-$independent constant $C>0$ such that $A\leq CB$ and $A\sim B$ means $A\lesssim B$ and $B\lesssim A$. $A\gtrsim B$ can defined similarly.

For function spaces,~$L^p(\mathbb{R}_+)(1\leq p\leq \infty)$~denotes
the usual Lebesgue space on~$\mathbb{R}_+$~with norm~$\|{\cdot}\|_{L^p}$
and for $k\in \mathbb{Z}_+$,~$H^k(\mathbb{R}_+)$~ represents the usual Sobolev space with the standard norm~$\|\cdot\|_k$. It is easy to see that $\|\cdot\|_0=\|\cdot\|_{L^2}$ and to simplify the notation, we set
~$\|\cdot\|:=\|\cdot\|_0=\|\cdot\|_{L^2}$ in the rest of this paper.

Finally, We denote by $C^k(I; H^p)$ the space of $k$-times continuously
differentiable functions on the interval $I$ with values in
$H^p(\mathbb{R}_+)$ and~$L^2(I; H^p)$ the space of~$L^2$-functions
on $I$ with values in~$H^p(\mathbb{R}_+)$.

\section{Preliminaries and main result}

To make the presentation easy to read, we divide this section into several subsections and the first one is on the properties of the viscous shock wave.

\subsection{Properties of the viscous shock wave}
In this subsection, we first recall some properties of the viscous shock wave. As already pointed out in the first section, we will only consider the $2-$viscous shock wave $(V_2^{\textrm{VSW}}(x-s_2t;w_-,w_+), U_2^{\textrm{VSW}}(x-st;w_-,w_+))$ connecting the states $w_-=(v_-,u_-)$ and $w_+=(v_+,u_+)$. To simplify the notations, we will set
$$
(V(\xi), U(\xi))\equiv\left(V_2^{\textrm{VSW}}(\xi;w_-,w_+), U_2^{\textrm{VSW}}(\xi;w_-,w_+)\right)
$$
and use $s$ to denote $s_2(v_-,v_+)$ in the rest of this paper. Consequently, we can deduce from \eqref{Viscous-shock-wave}, \eqref{wave speed}, \eqref{entropy condition} that $(V(\xi), U(\xi))$ solves
\begin{equation} \label{Viscous-shock-wave-2}
  \left\{
  \begin{array}{rl}
   U(\xi)-u_\pm&=-s\left(V(\xi)-v_\pm\right), \\[2mm]
   \frac{s\mu}{V(\xi)}\frac{d V(\xi)}{d\xi}&=-s^2\left(V
   (\xi)-v_\pm\right)
   -\left(p\left(V(\xi)\right)-p(v_\pm)\right)\\[2mm]
   &:=h\left(V(\xi)\right),\\[2mm]
   V(-\infty)&=v_-,\quad V(+\infty)=v_+,
  \end{array}
  \right.
\end{equation}
where $(v_-,u_-)$, $(v_+,u_+)$ and $s$ are assumed to satisfied the Rankine-Hugoniot condition
\begin{equation}\label{R-H}
\left\{
\begin{array}{rl}
s(v_+-v_-)&=u_--u_+,\\[2mm]
s(u_+-u_-)&=p(v_+)-p(v_-)
\end{array}
\right.
\end{equation}
and the entropy condition
\begin{equation}\label{entropy condition-2}
u_+<u_-.
\end{equation}

Recall that $\delta:=|v_+-v_-|$ denotes the strength of the $2-$viscous shock wave $(V(\xi), U(\xi))$, we have the following result on the $2-$viscous shock wave $(V(\xi), U(\xi))$:
\begin{Proposition}(cf. \cite{Huang-Matsumura-Shi-CMP-2003})
For any $(v_+, u_+), (v_-, u_-), s>0$ satisfying $v_+>v_->0$, the
Rankine-Hugoniot condition \eqref{R-H}, and the entropy condition \eqref{entropy condition-2}, there exists  a unique viscous shock wave $(V(\xi,U(\xi)) (\xi=x-st)$ up to a shift, which connects $(v_+,u_+)$, $(v_-,u_-)$ and satisfies
\begin{equation} \label{2.4}
 \begin{split}
&0<v_-<V(\xi)<v_+,\quad  u_+<U(\xi)<u_-,\\
&h(V(\xi))>0,\quad \frac{dV(\xi)}{d\xi}= \frac{V(\xi)h(V(\xi))}{s\mu} >0,\\
&\left|\Big(V(\xi)-v_\pm, U(\xi)-u_\pm\Big)\right|\leq O(1)\delta e^{-c_\pm |\xi|},\\
&\left|\left(\frac{dV(\xi)}{d\xi}, \frac{dU(\xi)}{d\xi}, \frac{d^2V(\xi)}{d\xi^2}, \frac{d^2U(\xi)}{d\xi^2}\right)\right|\leq O(1)\delta^2 e^{-c_\pm |\xi|}.
\end{split}
\end{equation}
Here $c_\pm=\frac{v_\pm |p'(v_\pm)+s^2|}{s\mu}>0$.

Moreover, if $v_-$ and $v_+$ are assumed to be independent of $\delta$, one can further deduce that the positive constant $O(1)$ in \eqref{2.4} depends only on $v_\pm$ but independent of $\delta$ and $c_\pm=O(1)\delta$ for some $\delta-$independent positive constant $O(1)$ which depends only on $v_\pm$.
\end{Proposition}

The proof of the above proposition is almost the same as the one given in \cite{Huang-Matsumura-Shi-CMP-2003}, what we want to emphasize here is that although we ask $\delta$, the strength of the viscous shock wave $(V(\xi), U(\xi))$, to be small, since $v_\pm$ are assumed to be independent of $\delta$, one can easily verify that the positive constant $O(1)$ appearing in \eqref{2.4} depends only on $v_\pm$ but is independent of $\delta$. Moreover, it is easy to verify that $c_\pm=O(1)\delta$ for some $\delta-$independent positive constant $O(1)$ by the Taylor formula. We thus omit the details for brevity.

\subsection{Main result}
Similar to that of \cite{Huang-Matsumura-Shi-CMP-2003}, we now make the following transformation
\begin{equation*}
 t=t, \quad \xi=x-s_-t
\end{equation*}
to transform the original problem \eqref{LagrangeI} to the following initial-boundary value problem
\begin{equation} \label{2.5}
  \begin{cases}
   v_t-s_{-}v_{\xi}-u_{\xi}=0, \quad  \xi\geq 0,\ t\geq 0,\\
   u_t-s_{-}u_{\xi}+p(v)_{\xi}=\mu\left(\frac{u_{\xi}}{v}\right)_{\xi}, \quad  \xi\geq 0,\ t\geq 0,\\
   (v(t,\xi),u(t,\xi))|_{\xi=0}=(v_-,u_-),\quad t\geq 0,\\
   (v(t,\xi),u(t,\xi))|_{t=0}=(v_0(\xi), u_0(\xi))\rightarrow(v_+,u_+),\ \textrm{as}\ 0\leq\xi\rightarrow +\infty.
  \end{cases}
\end{equation}

Since we will only consider the case $(v_-,u_-)\in \Omega_{\textrm{sub}}$, $(v_+,u_+) \in \textrm{S}_2(v_-,u_-)$, as classified in \cite{Matsumura-MAA-2001}, the large time behavior of the solution to (\ref{2.5}) is expected to be the suitably shifted $2-$viscous
shock wave $(V(\xi-(s-s_-)t+\sigma-\beta),U(\xi-(s-s_-)t+\sigma)-\beta)$ for some suitably chosen constant $\beta>0$. Here the shift $\sigma$ is defined as in \cite{Matsumura-Nishihara-CMP-1992}:
\begin{equation} \label{2.6}
\sigma=\frac{1}{v_+-v_-}\left\{\int^\infty_0\Big(v_0(\xi)-V(\xi-\beta)\Big)d\xi
-(s-s_-)\int^\infty_0\Big(V((s_--s)\tau-\beta)-v_-\Big)d\tau\right\},
\end{equation}
where $\beta>0$ is a suitably chosen constant whose precise range will be specified later. What we want to emphasize here is that the introduction of such a parameter $\beta$ is motivated by \cite{Matsumura-Mei-ARMA-1999} and the main purpose is to use it to control the boundary terms induced by the inflow boundary condition \eqref{2.5}$_3$ suitably.

By choosing the shift $\sigma$ as above, we can put
\begin{equation}\label{2.7}
  \begin{split}
  & \phi(t,\xi)=-\int^\infty_\xi[v(t,y)-V(y-(s-s_-)t+\sigma-\beta)]dy,\\
  & \psi(t,\xi)=-\int^\infty_\xi[u(t,y)-U(y-(s-s_-)t+\sigma-\beta)]dy,
  \end{split}
\end{equation}
which means
\begin{equation}\label{2.8}
  \begin{split}
  & v(t,\xi)=\phi_\xi(t,\xi)+V(\xi-(s-s_-)t+\sigma-\beta),\\
  & u(t,\xi)=\psi_\xi(t,\xi)+U(\xi-(s-s_-)t+\sigma-\beta),
  \end{split}
\end{equation}
and consequently the problem \eqref{2.5} can be reformulated as
\begin{equation} \label{2.9}
  \begin{cases}
  \phi_t-s_{-}\phi_{\xi}-\psi_{\xi}=0, \quad  \xi>0,\ t>0,\\
  \psi_t-s_{-}\psi_{\xi}+p(V+\phi_\xi)-p(V)
  =\mu\left(\frac{U_{\xi}+\psi_{\xi\xi}}{V+\phi_\xi}-\frac{U_{\xi}}{V}\right), \quad  \xi>0,\ t>0,\\
(\phi(t,\xi),\psi(t,\xi))|_{t=0}=(\phi_0(\xi),\psi_0(\xi))\rightarrow(0,0), \quad \textrm{as}\ \ 0\leq\xi\rightarrow +\infty,\\
  \phi(t,\xi)|_{\xi=0}=A(t):=-(s-s_-){\displaystyle\int^\infty_t}[V(-(s-s_-)\tau+\sigma-\beta)-v_-]d\tau,\quad t\geq 0.
  \end{cases}
\end{equation}

Now we turn to state our main result. Firstly, assume that
\begin{itemize}
\item [(H$_1$)] There exist $\delta-$independent constants $l\geq0$ and $C_0>0$ such that
\begin{equation}\label{2.10}
C^{-1}_0 \delta^ l \leq v_0(x) \leq C_0 (1+\delta^ {-l});
\end{equation}
\item [(H$_2$)] $(v_-,u_-)\in \Omega_{\textrm{sub}}$, $(v_+,u_+) \in \textrm{S}_2(v_-,u_-)$, $v_-$ and $v_+$ are positive constants independent of $\delta$;
\item[(H$_3$)] The initial data $(v_0(x),u_0(x))$ are assumed to satisfy
\begin{equation}\label{2.11}
\begin{split}
&\big(v_0(\xi)-V(\xi+\sigma-\beta),  u_0(\xi)-U(\xi+\sigma-\beta)\big)\in H^2(\mathbb{R}_+)\cap L^1 (\mathbb{R}_+),\\
& (\phi_0(\xi),\psi_0(\xi))\in  L^2 (\mathbb{R}_+).
\end{split}
\end{equation}
and the compatibility condition
\begin{equation}\label{2.12}
v_0(0)=v_-, \quad  u_0(0)=u_-.
\end{equation}
\end{itemize}

Under the above assumptions, we have
\begin{Theorem}
For any $(v_-,u_-)\in \Omega_{\textrm{sub}}$, $(v_+,u_+) \in \textrm{S}_2(v_-,u_-)$, $\gamma>1$, in addition to the assumptions
(H$_1$)-(H$_3$), we assume further that
\begin{equation}\label{2.13}
\begin{split}
& \left\|(\phi_0,\psi_0)\right\|_1\lesssim \delta^\alpha,\\
& \left\|\Big(\phi_{0xx},\psi_{0xx}\Big)\right\|\lesssim \left(1+\delta^{-\kappa}\right),\\
& |s_-|\sim u_-\sim \delta^h,\quad \beta=o(\delta^{-1})
\end{split}
\end{equation}
hold for some $\delta-$independent constants $C_1>0$, $\alpha>0, h>0$ and $\kappa>0$.
If the positive parameters $l$, $\alpha, h$ and $\kappa$ are assumed to satisfy:
\begin{equation}\label{2.14}
\begin{cases}
(\gamma+2)l<1,\quad (6\gamma+4)l<\alpha<\kappa,\\[2mm]
\frac{(\gamma+2)}{2}l<h <  \frac{7}{4}\big(\alpha-\frac{(\gamma+1)l}{2}\big),\\[2mm]
0<\theta < \min\left\{\frac{\gamma-1}{4(\gamma^2+3\gamma-2)}\left(\alpha
-\frac{(\gamma+1)l}{2}\right)
,\frac{\gamma-1}{\gamma^2+\gamma+2},
\frac{\gamma-1}{\gamma^2}h\right\},
\end{cases}
\end{equation}
where $\theta:=\kappa+l-\big(\alpha-\frac{(\gamma+1)l}{2}\big)$, then there exists a suitably small $\delta_0$ such that if $0<\delta\leq \delta_0$, the initial-boundary value problem \eqref{LagrangeI} has a unique solution $(v(t,x), u(t,x))$ satisfying
\begin{equation}\label{2.15}
\begin{split}
&\big(v(t,x)-V(x-st+\sigma-\beta), u(t,x)-U(x-st+\sigma-\beta)\big)\in C([0,+\infty), H^2(\mathbb{R}_+)),\\
&v_x(t,x)-V'(x-st+\sigma-\beta)\in  L^2 \big([0,+\infty), H^1(\mathbb{R}_+)\big),\\
&u_x(t,x)-U'(x-st+\sigma-\beta)\in  L^2 \big([0,+\infty), H^2(\mathbb{R}_+)\big),
\end{split}
\end{equation}
and
\begin{equation}\label{2.16}
C_2^{-1}\delta^{\frac{2\theta}{\gamma-1}}
\leq  v(t,x) \leq C_2\delta^{ -2\theta}
\end{equation}
for some positive constant $C_2$ independent of $\delta$. Furthermore, it holds that
\begin{equation} \label{largetime}
 \lim\limits_{t\to\infty} \sup_{x\geq s_-t}\left\{\Big|(v(t,x)-V(x-st+\sigma-\beta),u(t,x)-U(x-st+\sigma-\beta))\Big|\right\}=0.
\end{equation}
\end{Theorem}

\begin{Remark} Several remarks concerning Theorem 1.1 are listed below:
\begin{itemize}
\item We affirm that the set of the parameters $\alpha>0$, $\kappa>0$, $l\geq 0$, and $h>0$ satisfying the above conditions is not empty. In fact, let $l=0$, \eqref{2.14} is equivalent to
\begin{equation}\label{2.18}
\begin{cases}
0<h < \frac{7}{4}\alpha,\\
\alpha< \kappa < \min\left\{\frac{\gamma-1}{4(\gamma^2+3\gamma-2)}\alpha
,\frac{\gamma-1}{\gamma^2+\gamma+2},
\frac{\gamma-1}{\gamma^2}h\right\}+\alpha,
\end{cases}
\end{equation}
thus, the existence of $\alpha,\kappa,h$ is easy to verify.
\item It is easy to construct some initial perturbation $(\phi_0(\xi),\psi_0(\xi))$ satisfying the conditions listed in Theorem 1. In fact for each function $(f(\xi), g(\xi)\in H^2(\mathbb{R}_+)$ with
    $$
    \left\{
    \begin{array}{rl}
    {\textrm{Osc}}\ f'(\xi)&:=\sup\limits_{\xi\in\mathbb{R}_+}\big\{f'(\xi)\big\}-
    \inf\limits_{\xi\in\mathbb{R}_+}\big\{f'(\xi)\big\}\equiv A_1>0,\\[1mm]
    {\textrm{Osc}}\ g'(\xi)&:=\sup\limits_{\xi\in\mathbb{R}_+}\big\{g'(\xi)\big\}-
    \inf\limits_{\xi\in\mathbb{R}_+}\big\{g'(\xi)\big\}\equiv A_2>0,
    \end{array}
    \right.
    $$
    and for each $\alpha, \beta$ satisfying the conditions listed in Theorem 1, if we set
$$
\phi_0(\xi)=\delta^{\frac{3\alpha+\kappa}{2}}f\left(\delta^{-\kappa-\alpha}\xi\right),\quad
\psi_0(\xi)=\delta^{\frac{3\alpha+\kappa}{2}}g\left(\delta^{-\kappa-\alpha}\xi\right),
$$
one can verify that such a $(\phi_0(\xi),\psi_0(\xi))$ satisfies all the conditions listed in Theorem 1.

For such chosen $(\phi_0(\xi),\psi_0(\xi))$, we can construct the initial data $(v_0(\xi), u_0(\xi))$ through \eqref{2.18} and noticing that
$$
\phi'_0(\xi)=\delta^{\frac{\alpha-\kappa}{2}}f'\left(\delta^{-\alpha-\kappa}\xi\right),\quad
\psi'_0(\xi)=\delta^{\frac{\alpha-\kappa}{2}}g'\left(\delta^{-\alpha-\kappa}\xi\right),
$$
we can get that
$$
{\textrm{Osc}}\ \phi_0'(\xi)=\delta^{\frac{\alpha-\kappa}{2}}A_1,\quad
{\textrm{Osc}}\ \psi'_0(\xi)=\delta^{\frac{\alpha-\kappa}{2}}A_2.
$$
Thus from \eqref{2.18} and the fact that ${\textrm{Osc}}\ V(\xi)=\delta, {\textrm{Osc}}\ U(\xi)=s\delta$ which are assumed to be sufficiently small, one can deduce that
$$
{\textrm{Osc}}\ v_0(\xi)\sim \delta^{\frac{\alpha-\kappa}{2}}A_1,\quad
{\textrm{Osc}}\ u_0(\xi)\sim\delta^{\frac{\alpha-\kappa}{2}}A_2,
$$
which are sufficiently large for small $\delta$ since the parameters $\alpha$ and $\kappa$ satisfies $0<\alpha<\kappa$.
Consequently, the assumptions we imposed on the initial perturbations in Theorem 1 can indeed allow the oscillations of both the initial specific volume $v_0(\xi)$ and the initial velocity $u_0(\xi)$ to be large. Moreover, from the estimate \eqref{2.16} and the facts that
$$
\frac{2}{1-\gamma} \left\{\left(\alpha-\frac{(\gamma+1)l}{2}\right)-(\kappa+l)\right\}>0, \quad 2\left\{\left(\alpha-\frac{(\gamma+1)l}{2}\right)-(\kappa+l)\right\}<0,
$$
one can easily deduce that for each $t>0$, ${\textrm{Osc}}\ v(t,\xi)$ can also be large for small $\delta$.
\item Unlike that of \cite{Huang-Matsumura-Shi-CMP-2003}, we use the smallness of both $\delta$ and $u_-$ to control the possible growth of the solutions to the inflow problem \eqref{LagrangeI} caused by both the nonlinearity of the equations \eqref{LagrangeI}$_1$-\eqref{LagrangeI}$_2$ and the inflow boundary condition \eqref{LagrangeI}$_3$. It is worth to pointing out that since our result holds for any $\gamma>1$, it seems natural to ask $\delta$, the strength of the viscous shock wave, to be small because for large $\gamma$, the condition \eqref{strength-VSW} imposed in \cite{Huang-Matsumura-Shi-CMP-2003} implies that $\delta$ is sufficiently small. It would be of some interest to see whether similar stability result holds or not for a class of large initial perturbations satisfying similar conditions imposed in Theorem 1 but the strength of the viscous shock wave is only assumed to satisfy the condition \eqref{strength-VSW} imposed in \cite{Huang-Matsumura-Shi-CMP-2003}.
\item For the inflow problem to the full compressible Navier-Stokes equations, some nonlinear stability results with small initial perturbation are obtained in \cite{Huang-Li-Shi-CMS-2010, Nakamura-Nishibata-JHDE-2011, Qin-Nonlinearity-2011, Qin-Wang-SIMA-2009, Qin-Wang-SIMA-2011}, we are convinced that some results similar to that of \cite{Fan-Liu-Wang-Zhao-JDE-2014} and this paper can also be obtained which can allow the initial density to have large oscillations.
\end{itemize}
\end{Remark}

\subsection{Main ideas used to deduce our main result}
To yield a global solvability result to the initial-boundary value problem \eqref{2.9}, in addition to the difficulty on the possible growth of its solution $(\phi(t,\xi), \psi(t,\xi))$ induced by the nonlinearity of the equations \eqref{2.9}$_1$-\eqref{2.9}$_2$, the another is how to control the boundary terms caused by the inflow boundary condition \eqref{2.9}$_3$. To outline our main ideas used in this manuscript, we first recall the main ideas used in \cite{Huang-Matsumura-Shi-CMP-2003} to overcome these two difficulties as in the following:
\begin{itemize}
\item The first is to use the smallness of $N(T) :=\sup\limits_{0\leq t\leq T}\|(\phi,\psi)(t)\|_{L^\infty}$ to deal with the possible growth of $(\phi(t,\xi), \psi(t,\xi))$ caused by the nonlinearity of the equations \eqref{2.9}$_1$-\eqref{2.9}$_2$. One of the key points in such an argument is that, based on the a priori assumption that $N(T)$ is sufficiently small, one can deduce a uniform lower and upper positive bounds on the specific volume $v(t,\xi)$ which plays an essential role in deducing the desired a priori estimates on $(\phi(t,\xi), \psi(t,\xi))$;
\item As for the control of the boundary terms, there are two main ingredients in the analysis of \cite{Huang-Matsumura-Shi-CMP-2003}: The first is to introduce a parameter $\beta$ in the shift $\sigma$ given by \eqref{2.6}. It is motivated by \cite{Matsumura-Mei-ARMA-1999} and the main observation is that if one chooses $\beta>0$ sufficiently large such that $\beta\gg|\sigma|$, which is guaranteed by the smallness of the initial perturbation imposed in \cite{Huang-Matsumura-Shi-CMP-2003}, then some terms arising from the boundary terms can be controlled suitably, cf. Lemma 4.1 in \cite{Huang-Matsumura-Shi-CMP-2003}. The another is the introduction of a new variable $\overline{\psi}= \psi-s_-\phi$ to control the term $\psi^2(t,0)$ which was one of the main difficulties in the study of the nonlinear stability of viscous shock wave for the inflow problem \eqref{LagrangeI} pointed out in \cite{Matsumura-Nishihara-CMP-2001}. The observation in \cite{Huang-Matsumura-Shi-CMP-2003} is that when the energy method is applied to the new reformulated system, the first energy inequality does not contain the term $\overline{\psi}^2(t,0)$ provided that $|s_-|$ is small enough, which means that the estimates for the term $\psi^2(t,0)$ could be exactly bypassed and thus the desired a priori estimates can be obtained.
\end{itemize}
Based on the above two points, the authors obtained the desired $H^2(\mathbb{R}_+)-$norm a priori energy type estimates on $(\phi(t,\xi), \psi(t,\xi))$ in \cite{Huang-Matsumura-Shi-CMP-2003} in terms of the initial perturbation $(\phi_0(\xi), \psi_0(\xi))$ and the factor $e^{-c_-\beta}$ which can be chosen as small as we wanted if one chooses $\beta>0$ suitably large. From which the corresponding nonlinear stability result with small initial perturbation can be obtained by the standard continuation argument.

But to deduce a global existence result to the initial-boundary value problem \eqref{2.9} with large data which can allow the initial density to have large oscillations, the argument used in \cite{Huang-Matsumura-Shi-CMP-2003} cannot be used any longer. Our main ideas to yield the desired nonlinear stability results are the following:
\begin{itemize}
\item The first is on the way to use the parameter $\beta$ to control certain boundary terms. Our main observation is that, even though  $\delta$, the strength of the viscous shock wave, is assumed to be small in our analysis, under the assumptions imposed in Theorem 1 on the initial perturbation $(\phi_0(\xi), \psi_0(\xi))$, we can indeed prove that the shift $\sigma$ given by \eqref{2.6} satisfies $|\sigma|\lesssim\delta^{-1}$ and consequently, if we choose $\beta\sim o(\delta^{-1})$, then the boundary terms can also be controlled suitably, cf. Lemma 3.2;
\item Unlike the analysis in \cite{Huang-Matsumura-Shi-CMP-2003} where the smallness of the $H^2(\mathbb{R}_+)-$norm of the initial perturbation is used to the possible growth of $(\phi(t,\xi), \psi(t,\xi))$ caused by the nonlinearity of the equations \eqref{2.9}, we use the smallness of $\delta$, the strength of the viscous shock wave, instead. The main difficulty lies in how to yield the uniform positive lower and upper bounds on the specific volume $v(t,\xi)$. It is worth to emphasize that Kanel's argument \cite{Kanel-Diff Eqns-1968} plays an important role in this step and it was to guarantee that the whole analysis to be carried out smoothly that we need to ask the parameters  $l, \alpha, h,$ and $\kappa$ to satisfy the conditions \eqref{2.14} in Theorem 1.
\end{itemize}

\section{The proof of our main result}
This section is devoted to proving our main result. To this end, for some positive constants $0<T\leq +\infty$, $m$ and $M$, we first give the set of functions $X_{m,M}(0,T)$ from which the solution to the initial-boundary value problem \eqref{2.9} is sought as follows:
\begin{equation*}
    X_{m,M}(0,T)
    =\left\{(\phi(t,\xi),\psi(t,\xi))\ \left|\
    \begin{array}{c}
    (\phi,\psi)\in C([0,T];H^2);\ {\phi}_{\xi}\in L^2(0,T;H^1),{\psi}_{\xi}\in L^2(0,T;H^2),\\
\sup\limits_{[0,T]\times\mathbb{R}_+}\Big\{(V+\phi)(t,\xi)\Big\}\leq M,\inf\limits_{[0,T]\times\mathbb{R}_+}\Big\{(V+\phi)(t,\xi)\Big\}\geq m
   \end{array}\right.
   \right\}.
\end{equation*}

Our main result will be proved by combining the local existence result to the initial-boundary value problem \eqref{2.9} with some a priori estimates. For this purpose, we first consider the local existence result in the following subsection.

\subsection{Local solvability result in $X_{m,M}(0,T)$}
The local-in-time existence of the solution $(\phi(t,\xi),\psi(t,\xi))$ to the initial-boundary value problem \eqref{2.9} in $X_{m,M}(0,T)$ has been studied in \cite{Huang-Matsumura-Shi-CMP-2003}, we thus cite the result there as in the following proposition:
\begin{Proposition}\label{locals}
Let $({\phi}_0(\xi),{\psi}_0(\xi))\in H^2(\mathbb{R}_{+})$.
If
$$
\sup\limits_{\xi\in{\mathbb{R}_+}}\big\{V(\xi)+\phi_0(\xi)\big\}\leq M,\quad
\inf\limits_{{\xi\in\mathbb{R}_+}}\big\{V(\xi)+\phi_0(\xi)\big\}\geq m,
$$
then there exists $t_0 >0$ depending only on $m$, $M$ and $\|({\phi}_0,{\psi}_{0})\|_2$ such that the initial-boundary value problem \eqref{2.9} has a unique solution $(\phi(t,\xi),\psi(t,\xi))\in X_{m/2,2M}(0,t_0)$ satisfying
\begin{equation} \label{3.1}
\begin{split}
\|(\phi,{\psi})(t)\|_1 \leq 2\|(\phi_0,\psi_0)\|_1,\\
 \|(\phi,{\psi})(t)\|_2 \leq 2\|(\phi_0,\psi_0)\|_2.
 \end{split}
\end{equation}
\end{Proposition}

\subsection{Certain a priori estimates}
Assume that the local solution $(\phi(t,\xi),{\psi}(t,\xi))$ obtained in the Proposition \ref{locals} has been extended to the time $t=T\geq t_0 \geq 0$, in order to show that $T=\infty$,  we now turn to deduce
certain energy type estimates on $(\phi(t,\xi),\psi(t,\xi))\in X_{m,M}(0,T)$ for some positive constants $m$ and $M$
and consequently
\begin{equation}\label{3.2}
\frac{1}{m}\leq V(\xi-(s-s_-)t+\sigma-\beta)+\phi_x(t,\xi) \leq M, \quad \forall(t,\xi)\in[0,T]\times \mathbb{R}_+.
\end{equation}
Without loss of generality, we can assume that $m\geq1$ and $M\geq1$.

Moreover we assume further that
\begin{equation}\label{L-infty-assumption}
\big|(\phi(t,\xi),\psi(t,\xi))\big|\leq N(T),\quad \forall (t,\xi)\in[0,T]\times\mathbb{R}_+.
\end{equation}

Firstly, we give some estimates on the shift $\sigma$ and the boundary terms.
\begin{Lemma} Assume that the conditions listed in Theorem in 2.1, it holds that $|\sigma| \lesssim \delta^{-1}$.
\end{Lemma}
\noindent{\bf Proof}: Noticing that $\frac \delta{c_\pm}\sim 1$, it is easy to see from Proposition 2.1 that
\begin{equation*}\label{3.0.1}
\int^\infty_0|v_+-V(y+d)|dy \lesssim \delta \int^\infty_0e^{-c_+|y+d|}dy  \lesssim \frac\delta{c_+}\lesssim 1
\end{equation*}
holds for any $d \in  \mathbb{R}$, we thus get that
\begin{equation} \label{3.0.2}
\begin{split}
|\sigma|
\lesssim & \frac{1}{\delta}\left\{\Bigg|\int^\infty_0\Big(v_0(\xi)-V(\xi+\sigma-\beta)\Big)d\xi\Bigg|
+\Bigg|\int^\infty_0\Big(V(\xi+\sigma-\beta)-V(\xi-\beta)\Big)d\xi\Bigg|\right\}\\
&+\frac{1}{\delta}\Bigg|-(s-s_-)\int^\infty_0\Big(V(-(s-s_-)\tau-\beta)-v_-\Big)d\tau\Bigg|\\
\lesssim & \frac{1}{\delta}\left\{|\phi_0(0)|+\Bigg|\int^\infty_0\Big(V(\xi+\sigma-\beta)-v_+\Big)d\xi\Bigg|
+\Bigg|\int^\infty_0\Big(v_+-V(\xi-\beta)\Big)d\xi\Bigg|\right\}\\
&+ \frac{1}{\delta}\Bigg|{\delta}\int^\infty_0e^{-c_-|-(s-s_-)\tau-\beta|}d(-(s-s_-)\tau)\Bigg|\\
\lesssim & \frac{1}{\delta}\left\{|\phi_0(0)|+1\right\}.
\end{split}
\end{equation}

Under the conditions listed in Theorem 1, it can be deduce from \eqref{3.0.2} that $|\sigma| \lesssim \delta^{-1}$.
This completes the proof of Lemma 3.1.\\

With the estimate on the shift $\sigma$ obtained in Lemma 3.1 in hand, we now turn to control the involving boundary terms suitably which is the main content of the next lemma
\begin{Lemma} Assume that $(\phi(t,\xi),{\psi}(t,\xi))$ obtained in the Proposition \ref{locals} has been extend to the time step $t=T$ and satisfies the conditions listed above, i.e. $(\phi(t,\xi),{\psi}(t,\xi))\in X_{m,M}(0,T)$ and the a priori assumption \eqref{L-infty-assumption} is assumed to hold, then if $\beta>0$ is chosen sufficiently large such that $\beta\geq|\sigma|$, it holds for each $0\leq t \leq T$ that
\begin{equation}\label{3.5}
\begin{split}
&\int^t_0|\phi(\tau,0)|d\tau \lesssim \delta^{-1}e^{-c_-\beta}, \quad \int^t_0|\phi_\xi(\tau,0)|d\tau \lesssim e^{-c_-\beta},\quad \int^t_0|\psi_\xi(\tau,0)|d\tau \lesssim e^{-c_-\beta}, \\
&\int^t_0|\phi_t(\tau,0)|d\tau \lesssim e^{-c_-\beta},\quad
\int^t_0|\phi_{t \xi}(\tau,0)|d\tau \lesssim \delta e^{-c_-\beta},\quad
\int^t_0|\psi_{t \xi}(\tau,0)|d\tau \lesssim \delta e^{-c_-\beta}.
\end{split}
\end{equation}
Although $c_-\sim\delta$, if we choose $\beta>0$ sufficiently large such that
$\delta\cdot\beta=o\left(\delta^{-1}\right)$, then if $\delta$, the strength of the viscous shock wave, is small enough, the estimates \eqref{3.5} tell us that the contributions of all these integrals involving the boundary terms can be as small as wanted.
\end{Lemma}
\noindent{\bf Proof}: Firstly, noticing that $c_-=O(1)\delta$ and $\beta\geq|\sigma|$, one can get from Proposition 2.1 that
\begin{eqnarray*}
|\phi(t,0)|=|A(t)|&=&(s-s_-)
\int^\infty_t\big|V(-(s-s_-)\tau+\sigma-\beta)-v_-\big|d\tau\\
&\lesssim &(s-s_-)\delta \int^\infty_te^{-c_-((s-s_-)\tau+\beta-\sigma)}d\tau\\
&\lesssim& \frac\delta{c_-} e^{-c_-((s-s_-)t-\sigma+\beta)}\\
&\lesssim& e^{-c_-((s-s_-)t-\sigma+\beta)}.
\end{eqnarray*}
Consequently, it holds
\begin{equation*}
\int^t_0|\phi(\tau,0)|d\tau\lesssim \int^t_0 e^{-c_-((s-s_-)\tau-\sigma+\beta)} d\tau
 \lesssim  c_-^{-1}e^{c_-(\sigma-\beta)} \lesssim  \delta^{-1}e^{-c_-\beta},
\end{equation*}
where we have used $|\sigma| \lesssim \delta^{-1}$.
On the other hand, due to
\begin{equation*}
\begin{split}
&\phi_\xi(t,0)= v(t,0)-V(-(s-s_-)t+\sigma-\beta)=v_--V(-(s-s_-)t+\sigma-\beta),\\
&\psi_\xi(t,0)= u(t,0)-U(-(s-s_-)t+\sigma-\beta)=u_--U(-(s-s_-)t+\sigma-\beta),
\end{split}
\end{equation*}
we can get from Proposition 2.1 that
\begin{eqnarray*}
\int^t_0|\phi_\xi(\tau,0)| d\tau &=&\int^t_0\big|v_--V(-(s-s_-)\tau+\sigma-\beta)\big|d\tau\\
&\lesssim&\delta\int^t_0 e^{-c_-|-(s-s_-)\tau+\sigma-\beta|} d\tau\\
&\lesssim& e^{-c_-\beta},\\
\int^t_0|\phi_{t\xi}(\tau,0)| d\tau &=&(s-s_-)\int^t_0\big|V'(-(s-s_-)\tau+\sigma-\beta)\big|d\tau\\
&\lesssim&\delta^2\int^t_0 e^{-c_-|(s-s_-)\tau+\sigma-\beta|} d\tau\\ &\lesssim&\delta e^{-c_-\beta}.
\end{eqnarray*}
Here again we have used the facts that $s$ is independent of $\delta$, $c_-=O(1)\delta$, $\beta\geq |\sigma|$, and $|\sigma|\lesssim \delta^{-1}$.

Similarly, one can get that
\begin{eqnarray*}
\int^t_0|\psi_\xi(\tau,0)| d\tau &=&\int^t_0\big|u_--U(-(s-s_-)\tau+\sigma-\beta)\big|d\tau\\
&\lesssim&\delta\int^t_0 e^{-c_-|-(s-s_-)\tau+\sigma-\beta|} d\tau\\
&\lesssim& e^{-c_-\beta},\\
\int^t_0|\psi_{t\xi}(\tau,0)| d\tau &=&(s-s_-)\int^t_0\big|U'(-(s-s_-)\tau+\sigma-\beta)\big|d\tau\\
&\lesssim&\delta^2\int^t_0 e^{-c_-|-(s-s_-)\tau+\sigma-\beta|} d\tau\\ &\lesssim&\delta e^{-c_-\beta}.
\end{eqnarray*}

At last, by the equation
\begin{equation*}
\phi_t(\tau,0) = s_- \phi_\xi(\tau,0)+ \psi_\xi(\tau,0),
\end{equation*}
it yields
\begin{equation*}
\int^t_0|\phi_t(\tau,0)| d\tau = \left|s_- \int^t_0\phi_\xi(\tau,0)d\tau+ \int^t_0\psi_\xi(\tau,0)d\tau\right|
\lesssim e^{-c_-\beta}.
\end{equation*}
This completes the proof of Lemma 3.2.\\

The next result is concerned with the basic energy estimate, which is stated in the following
lemma.
\begin{Lemma} Let $(\phi(t,\xi),{\psi}(t,\xi)$ satisfy the conditions listed in Lemma 3.2, then there exists a sufficiently small positive constant $\epsilon_1$ independent of
$\delta$ such that if
$$N(t)<\epsilon_1, \quad |s_-|<\epsilon_1,
$$
it holds for each $0\leq t \leq T$,
\begin{equation}\label{3.6}
\begin{split}
&\|(\phi,\psi)(t)\|^2 + \int^t_0\int^\infty_{0}\psi^2_\xi d\xi d\tau \\
\lesssim &\|(\phi_0,\psi_0)\|^2+\delta^{-1}e^{-c_-\beta}+ N(t)^{\frac{2}{3}}\int^t_0\int_{\mathbb{R}_+}\frac{\psi^2_{\xi\xi}}{v} d\xi d\tau\\
&+\left(m^{\gamma+2}N(t)+m^2\delta+s^2_-\right) \int^t_0\int_{\mathbb{R}_+}\phi^2_{\xi} d\xi d\tau.
\end{split}
\end{equation}
\end{Lemma}
\noindent{\bf Proof}: As in \cite{Huang-Matsumura-Shi-CMP-2003}, let $ \overline{\psi}= \psi-s_-\phi$, the problem \eqref{2.9} can be changed into
\begin{equation} \label{3.7}
  \begin{cases}
  \phi_t-2s_{-}\phi_{\xi}-\overline{\psi}_{\xi}=0, \quad  \xi>0,\ t>0,\\
  \overline{\psi}_t-\left(-p'(V)-s^2_{-}\right)\phi_{\xi}-\mu\frac{\psi_{\xi\xi}}{V}=F, \quad  \xi>0,\ t>0,\\
\left(\phi(0,\xi),\overline{\psi}(0,\xi)\right)=(\phi_0(\xi),\psi_0(\xi)-s_-\phi_0(\xi)),\quad \xi\geq 0,\\
  \phi(t,0)=A(t),  \psi_{\xi}(t,0)=A'(t),\quad t\geq 0,
  \end{cases}
\end{equation}
where
\begin{equation*}
 F=-\big(p(V+\phi_\xi)-p(V)-p'(V)\phi_\xi\big)-\frac{\mu U_\xi \phi_\xi}{vV}-\frac{\mu \phi_\xi {\psi_{\xi\xi}}}{vV}.
\end{equation*}

Multiplying $(\ref{3.7})_1$ by $\phi$ and $(\ref{3.7})_2$ by $(-p'(V)-s^2_-)^{-1}\overline{\psi}$,
then we have
\begin{equation}\label{3.8}
\begin{split}
&\left(\frac{\phi^2}{2}-\frac{\overline{\psi}^2}{2(p'(V)+s^2_-)}\right)_t
+\left(\frac{1}{2(p'(V)+s^2_-)}\right)_t\overline{\psi}^2 -\frac{\mu \overline{\psi}_\xi\psi_\xi}{V(p'(V)+s^2_-)}\\
&+\left\{-s_-\phi^2-\phi\overline{\psi}-\frac{\mu \overline{\psi}\psi_\xi}{V(p'(V)+s^2_-)}\right\}_\xi \\
=& \left\{\frac{\mu}{V(p'(V)+s^2_-)}\right\}_\xi\overline{\psi}\psi_\xi
+ \frac{\mu U_\xi \phi_\xi \overline{\psi}}{vV(p'(V)+s^2_-)}
+\frac{\mu  \phi_\xi \overline{\psi}{\psi}_{\xi\xi}}{vV(p'(V)+s^2_-)}\\
&+\frac{\mu \overline{\psi}}{(p'(V)+s^2_-)}\big[p(V+\phi_\xi)-p(V)-p'(V)\phi_\xi\big].
\end{split}
\end{equation}

Since $|s_-|$ is small, $v_-<V(\xi)<v_+$ with $v_\pm$ independent of $\delta$, $ p'(V)<0,$ $p''(V)>0,$
$$\frac{\partial V(\xi-(s-s_-)t+\sigma-\beta)}{\partial t}=-(s-s_-)V'(\xi-(s-s_-)t+\sigma-\beta)<0,
$$
and noticing the fact that $\overline{\psi}=\psi-s_-\phi$ and $|p'(V)+s_-^2|$ can be bounded by some positive constant independent of $\delta$ from both below and above, which follows from the facts that both $\delta$ and $|s_-|$ are assumed to be sufficiently small and $v_\pm$, the far fields of $V(\xi)$, are independent of $\delta$,
then integrating the above identity with respect to
$t$ and $x$ over $[0,t]\times \mathbb{R}_+$ yields
\begin{equation}\label{3.9}
\begin{split}
&\int_{\mathbb{R}_+}\bigg(\frac{\phi^2}{2}-\frac{\overline{\psi}^2}
{2(p'(V)+s^2_-)}\bigg) d\xi
+ \int^t_0\int_{\mathbb{R}_+}\left(|V_t|\overline{\psi}^2+\frac{\mu \psi^2_\xi}{V(p'(V)+s^2_-)}\right) d\xi d\tau \\
\lesssim & \|(\phi_0,\overline{\psi}_0)\|^2+ \underbrace{\left|\int^t_0\left(s_-\phi^2+\phi\overline{\psi}+
\frac{\mu\overline{\psi}{\psi}_\xi}{V(s_-^2+p'(V))}\right)(\tau,0) d\tau\right|}_{I_1}\\
&+\underbrace{\left|\int^t_0\int_{\mathbb{R}_+}\frac{\mu  \phi_\xi \overline{\psi}{{\psi}}_{\xi\xi}}{vV(p'(V)+s^2_-)}
 d\xi d\tau\right|}_{I_2}+\int^t_0\int_{\mathbb{R}_+}\frac{\mu s^2_- \phi^2_\xi}{V(p'(V)+s^2_-)}
 d\xi d\tau\\
 &+\underbrace{\left|\int^t_0\int_{\mathbb{R}_+}\left(\overline{\psi}
 \left(p(V+\phi_\xi)-p(V)-p'(V)\phi_\xi\right)\right)
  d\xi d\tau\right|}_{I_3}\\
&+\underbrace{\left|\int^t_0\int_{\mathbb{R}_+}\left(
\left\{\frac{\mu}{V(p'(V)+s^2_-)}\right\}_\xi\overline{\psi}{\psi}_\xi
+ \frac{\mu U_\xi \phi_\xi \overline{\psi}}{vV(p'(V)+s^2_-)}\right) d\xi d\tau\right|}_{I_4}.
\end{split}
\end{equation}

Now we deal with $I_j (j=1,2,3,4)$ term by term. To this end, noticing first that $\|\overline{\psi}\|_{L^\infty}=\|\psi-s_-\phi\|_{L^\infty}\lesssim N(t)$, we can get from Lemma 3.2 that
\begin{equation}\label{3.10}
\begin{split}
I_1&\lesssim N(t)\int^t_0\left(|s_-\phi|+|\phi|+\left|\frac{\mu{\psi}_\xi}{V(s_-^2+p'(V))}\right|\right)(\tau,0) d\tau\\
&\lesssim N(t)\int^t_0\left(|\phi|+\left|\psi_\xi\right|\right)(\tau,0) d\tau \lesssim \delta^{-1}e^{-c_-\beta},\\
{I_2}&\lesssim N(t)^{\frac{2}{3}}\int^t_0\int^\infty_0\frac{{\psi}^2_{\xi\xi}}{v} d\xi d\tau
+N(t)^{\frac{1}{3}}m \int^t_0 \int^\infty_0\phi^2_\xi  d\xi d\tau.
\end{split}
\end{equation}
On the other hand, due to
\begin{equation}\label{3.11}
\begin{split}
\left|p(v)-p(V)-p'(V)\phi_\xi\right|
&=\phi^2_\xi\left|\int^1_0\int^1_0p''(\theta_1\theta_2v+(1-\theta_1\theta_2)V)d\theta_1 d\theta_2\right|\\
&\lesssim \left(v^{-\gamma-2}+V^{-\gamma-2}\right)\phi^2_\xi,
\end{split}
\end{equation}
one can deduce from \eqref{3.2} again that
\begin{equation*}
{I_3}\lesssim N(t)\int^t_0\int^\infty_0\left(v^{-\gamma-2}+V^{-\gamma-2}\right)\phi^2_\xi d \xi d\tau
\lesssim N(t)m^{\gamma+2}\int^t_0\int^\infty_0\phi^2_\xi d \xi d\tau
\end{equation*}
and
\begin{equation*}
{I_4}\lesssim \eta\int^t_0\int^\infty_0|V_t|\overline{\psi}^2 d \xi d\tau
+\delta\int^t_0\int^\infty_0 \frac{\mu \psi^2_\xi}{V(p'(V)+s^2_-)}d \xi d\tau+\delta m^2\int^t_0\int_{\mathbb{R}_+}\phi^2_\xi d \xi d\tau.
\end{equation*}
By choosing $\eta>0$ small enough, we can get by inserting the above estimates on $I_j (j=1,2,3,4)$ into \eqref{3.9} that
\begin{equation}\label{3.12}
\begin{split}
&\left\|\left(\phi,\overline{\psi}\right)(t)\right\|^2
+ \int^t_0\int^\infty_0\left(|V_t|\overline{\psi}^2+{\psi}^2_\xi \right)  d\xi d\tau \\
\lesssim &\|(\phi_0,\overline{\psi}_0)\|^2+\delta^{-1}e^{-c_-\beta}
+N(t)^{\frac{2}{3}}\int^t_0\int^\infty_0\frac{\psi^2_{\xi\xi}}{v} d\xi d\tau\\
&
+\left(N(t)m^{\gamma+2}+\delta m^2+s^2_-\right) \int^t_0\int^\infty_0\phi^2_{\xi} d\xi d\tau.
\end{split}
\end{equation}
Recall that $\overline{\psi}=\psi-s_-\phi$, it is easy to get the estimate \eqref{3.6} from \eqref{3.12}. This completes the proof of Lemma 3.3.\\

Now we turn to deduce the higher order energy estimates on $(\phi(t,\xi),\psi(t,\xi))$. To this end, we can get first that
\begin{Lemma} Under the same assumptions listed in Lemma 3.2,  it holds for each $0\leq t \leq T$ that
\begin{equation}\label{3.13}
\begin{split}
&\left\|\left(\sqrt{\Phi},\psi_\xi\right)(t)\right\|^2
+ \int^t_0\int^\infty_0 \frac{{\psi}^2_{\xi\xi}}{v}d\xi d\tau \\
\lesssim & \left\|\left(\sqrt{\Phi_0},\psi_{0\xi}\right)\right\|^2+ e^{-c_-\beta}
+\delta^2 m^{\gamma+2} \int^t_0\int^\infty_0 \phi^2_{\xi} d\xi d\tau.
\end{split}
\end{equation}
\end{Lemma}
\noindent{\bf Proof}: Differentiating (\ref{2.9})$_1$-\eqref{2.9}$_2$ with respect to $\xi$ once yields
\begin{equation} \label{3.14}
  \begin{cases}
  \phi_{\xi t}-s_{-}\phi_{\xi\xi}-\psi_{\xi\xi}=0, \quad  \xi>0,\ t>0,\\
  \psi_{\xi t}-s_{-}\psi_{\xi\xi}+(p(v)-p(V))_\xi-\mu\left(\frac{\psi_{\xi\xi}}{v}\right)_\xi
  =\mu\left(\frac{U_\xi}{v}-\frac{U_\xi}{V}\right)_\xi,\quad  \xi>0,\ t>0.
  \end{cases}
\end{equation}

Multiplying $(\ref{3.14})_1$ by $(p(V)-p(v))$ and $(\ref{3.14})_2$ by $\psi_\xi$, it holds that
\begin{equation}\label{3.15}
\begin{split}
&\left(\Phi+\frac{\psi^2_\xi}{2}\right)_t+\mu\frac{\psi^2_{\xi\xi}}{v}
+\left(-\frac{s_-}{2}\psi^2_\xi+(p(v)-p(V))\psi_\xi
-\mu\left(\frac{U_\xi+\psi_{\xi\xi}}{v}-\frac{U_\xi}{V}\right)\psi_\xi\right)_\xi\\
&=\mu\frac{U_\xi\phi_\xi\psi_{\xi\xi}}{vV}-V_t(p(v)-p(V)-p'(V)\phi_\xi),
\end{split}
\end{equation}
where
\begin{equation*}
\Phi(v,V)=p(V)(v-V)-\int^v_V p(\eta) d\eta.
\end{equation*}

Integrating the above identity with respect to $t$ and $x$ over $[0,t]\times \mathbb{R}_+$ yields
\begin{equation}\label{3.16}
\begin{split}
&\left\|\left(\sqrt{\Phi},\psi_\xi\right)(t)\right\|^2
+ \mu\int^t_0\int^\infty_0 \frac{{\psi}^2_{\xi\xi}}{v}d\xi d\tau \\
\lesssim & \left\|\left(\sqrt{\Phi_0},\psi_{0\xi}\right)\right\|^2
+ \underbrace{\left|\int^t_0\bigg(-\frac{s_-}{2}\psi^2_\xi+(p(v)-p(V))\psi_\xi
-\mu\left(\frac{U_\xi+\psi_{\xi\xi}}{v}-\frac{U_\xi}{V}\right)\psi_\xi\bigg)(\tau,0) d\tau\right|}_{I_5}\\
&+\underbrace{\left|\int^t_0\int^\infty_0\frac{U_\xi  \phi_\xi {\psi}_{\xi\xi}}{vV} d\xi d\tau\right|}_{I_6}
+\underbrace{\left|\int^t_0\int^\infty_0V_t\big(p(V+\phi_\xi)-p(V)-p'(V)\phi_\xi\big) d\xi d\tau\right|}_{I_7}.
\end{split}
\end{equation}

To bound $I_j (j=5,6,7)$ term by term, noticing first that
\begin{equation*}
 \mu\frac{\psi_{\xi\xi}}{v}= \psi_{t}-s_{-}\psi_{\xi}+(p(v)-p(V))+\mu\frac{U_\xi \phi_\xi}{vV},
\end{equation*}
we have
\begin{equation*}
\begin{split}
 \mu\frac{\psi_{\xi\xi}\psi_{\xi}}{v}
 = (\psi\psi_{\xi})_t-\psi\psi_{t\xi}-s_{-}\psi^2_{\xi}+(p(v)-p(V))\psi_{\xi}
 +\mu\frac{U_\xi \phi_\xi\psi_{\xi}}{vV},
 \end{split}
\end{equation*}
and consequently
\begin{equation*}
\begin{split}
 &\int^t_0\bigg(\mu\frac{\psi_{\xi\xi}\psi_{\xi}}{v}-(p(v)-p(V))\psi_\xi-\frac{\mu U_\xi \phi_\xi\psi_{\xi}}{vV}\bigg)(\tau,0)d\tau\\
 =& \psi\psi_{\xi}(t,0)-\psi\psi_{\xi}(0,0)- \int^t_0\psi\psi_{t\xi}(\tau,0)d\tau
 -s_{-} \int^t_0\psi^2_{\xi}(\tau,0)d\tau.
 \end{split}
\end{equation*}
Thus we can get from Lemma 3.2 that
\begin{equation}\label{3.17}
\begin{split}
I_5 &= \left|\frac{s_{-}}{2} \int^t_0\psi^2_{\xi}(\tau,0)d\tau-\psi(t,0)\psi_{\xi}(t,0)+\psi(0,0)\psi_{\xi}(0,0)+ \int^t_0\left(\psi\psi_{t\xi}\right)(\tau,0)d\tau\right|\\
&\lesssim (\delta + N(t)) e^{-c_-\beta} \lesssim e^{-c_-\beta},\\
I_6 & \lesssim \frac{\mu}{2} \int^t_0\int^\infty_0 \frac{{\psi}^2_{\xi\xi}}{v}d\xi d\tau
+\delta^4 m \int^t_0\int^\infty_0 \phi^2_{\xi}d\xi d\tau,\\
I_7 & \lesssim\delta^2 m^{\gamma+2} \int^t_0\int^\infty_0 \phi^2_{\xi}d\xi d\tau.
\end{split}
\end{equation}
Inserting the above estimates into  (\ref{3.16}), we can get (\ref{3.13}) immediately. This completes the proof of Lemma 3.4.\\

Combing the estimates obtained in Lemma 3.3 and Lemma 3.4, it holds
\begin{equation}\label{3.18}
\begin{split}
&\left\|\left(\phi,\psi,\sqrt{\Phi},\psi_\xi\right)(t)\right\|^2
+ \int^t_0\int^\infty_0\left({\psi}^2_{\xi}+ \frac{{\psi}^2_{\xi\xi}}{v}\right)d\xi d\tau \\
\lesssim & \left\|\left(\phi_0,\psi_0,\sqrt{\Phi_0},\psi_{0\xi}\right)\right\|^2+ \delta^{-1}e^{-c_-\beta}+\left(\left(N(t)+\delta^2\right) m^{\gamma+2} +s^2_-\right)\int^t_0\int^\infty_0 \phi^2_{\xi} d\xi d\tau.
\end{split}
\end{equation}

Now, we should deal with the term $\int^t_0\int^\infty_0\phi^2_{\xi}d\xi d\tau$. To this end, due to
\begin{equation*}
\psi_t-s_{-}\psi_{\xi}+(p(V+\phi_\xi)-p(V))-\mu\frac{\psi_{\xi\xi}}{v}
  =\mu\left(\frac{U_\xi}{v}-\frac{U_\xi}{V}\right),
\end{equation*}
we can get by multiplying the above equation by $\phi_\xi$ and by noticing $U_\xi<0$ that
\begin{equation}\label{3.19}
\begin{split}
&(\psi \phi_\xi)_t-\{\psi(s_{-}\phi_\xi+\psi_{\xi})\}_\xi+\psi^2_{\xi}
-\mu\frac{\psi_{\xi\xi}}{v} \phi_\xi\\
  =&-(p(V+\phi_\xi)-p(V))\phi_\xi-\mu\frac{U_\xi\phi^2_\xi}{v V}\\
   \geq& -(p(V+\phi_\xi)-p(V))\phi_\xi\\
  =&  -\left(\int^1_0p'(V+\theta\phi_\xi)d\theta\right)\phi^2_\xi\gtrsim M^{-\gamma-1} \phi^2_\xi.
 \end{split}
\end{equation}

Integrating the above inequality with respect to $t$
and $x$ over $[0,t]\times \mathbb{R}_+$, we get
\begin{equation}\label{3.20}
\begin{split}
M^{-\gamma-1} \int^t_0\int^\infty_{0} \phi^2_\xi d\xi d\tau
\lesssim& \|\psi\phi_\xi(t)\|+\|\psi_0\phi_{0\xi}\|
+\left|\int^t_0\bigg(\psi(s_{-}\phi_\xi+\psi_{\xi})\bigg)(\tau,0) d\tau\right| \\
&+ \int^t_0\int^\infty_{0}\psi^2_\xi d\xi d\tau
+m\int^t_0\int^\infty_{0} \frac{\psi^2_{\xi\xi}}{v} d\xi d\tau.
\end{split}
\end{equation}
Thus we finally get from \eqref{3.18}, \eqref{3.20}, and Lemma 3.2 that
\begin{equation}\label{3.21}
\begin{split}
&\int^t_0\int^\infty_{0} \phi^2_\xi d\xi d\tau\\
\lesssim & M^{\gamma+1}\Big(\|\psi(t)\|\|\phi_\xi(t)\|+\|\psi_0\|\|\phi_{0\xi}\|+e^{-c_-\beta}\Big)
+ mM^{\gamma+1}\int^t_0\int^\infty_0 \left( \psi^2_\xi +\frac{\psi^2_{\xi\xi}}{v} \right) d\xi d\tau\\
\lesssim & mM^{\gamma+1}\left(\left\| \left(\phi_0,\psi_0,\sqrt{\Phi_0},\psi_{0\xi}\right)\right\|^2+\delta^{-1} e^{-c_-\beta}
\right)\\
&+\left(\left(N(t)+\delta^2\right) m^{\gamma+3}M^{\gamma+1} +s^2_-mM^{\gamma+1}\right)\int^t_0\int^\infty_{0} \phi^2_{\xi}d\xi d\tau.
\end{split}
\end{equation}
If $N(t)$, $|s_-|$ and $\delta$ are chosen sufficiently small such that
\begin{equation*}
\left(N(t)+\delta^2\right) m^{\gamma+3}M^{\gamma+1}\leq \epsilon_2,\quad s^2_-mM^{\gamma+1} \leq \epsilon_2
\end{equation*}
hold for some sufficiently small positive constant $\epsilon_2>0$ independent $\delta$, then we can deduce from \eqref{3.21} that
\begin{equation*}
\begin{split}
\int^t_0\int^\infty_{0} \phi^2_\xi d\xi d\tau
\lesssim mM^{\gamma+1}\left(\left\|\left(\phi_0,\psi_0,\sqrt{\Phi_0},\psi_{0\xi}\right)\right\|^2+ \delta^{-1}e^{-c_-\beta}\right).
\end{split}
\end{equation*}

Inserting the above estimate into  (\ref{3.18}), we can get the following result
\begin{Lemma} Under the assumptions listed in Lemma 3.2, there exists a sufficiently small positive constant $\epsilon_2$ independent of $\delta$ such that if
\begin{equation}\label{3.22}
\left(N(t)+\delta^2\right) m^{\gamma+3}M^{\gamma+1}\leq \epsilon_2,\quad s^2_-mM^{\gamma+1} \leq \epsilon_2,
\end{equation}
then it holds for each $0\leq t \leq T$ that
\begin{equation}\label{3.23}
\begin{split}
\int^t_0\int^\infty_{0} \phi^2_\xi d\xi d\tau
\lesssim mM^{\gamma+1}\left(\left\|\left(\phi_0,\psi_0,\sqrt{\Phi_0},\psi_{0\xi}\right)\right\|^2+ \delta^{-1}e^{-c_-\beta}
\right)
\end{split}
\end{equation}
and
\begin{equation}\label{3.24}
\begin{split}
&\left\|\left(\phi,\psi,\sqrt{\Phi},\psi_\xi\right)(t)\right\|^2
+ \int^t_0\int^\infty_{0}\left({\psi}^2_{\xi}+ \frac{{\psi}^2_{\xi\xi}}{v}\right)d\xi d\tau \\
\lesssim & \left\|\left(\phi_0,\psi_0,\sqrt{\Phi_0},\psi_{0\xi}\right)\right\|^2
+\delta^{-1}e^{-c_-\beta}.
\end{split}
\end{equation}
\end{Lemma}

Now we pay attention to the term $\int^t_0\int^\infty_0 \frac{\phi^2_{\xi\xi}}{v}d\xi d\tau$. To this end, we set $\widetilde{v}:=\frac{v}{V}$,
then $\Phi(v,V)$ can be rewritten as
\begin{equation}\label{3.25}
\Phi(v,V)= V^{\gamma+1}\widetilde{\Phi}(\widetilde{v}),\quad
\widetilde{\Phi}(\widetilde{v})=\widetilde{v}-1+\frac{1}{\gamma-1}\left
(\widetilde{v}^{-\gamma+1}-1\right).
\end{equation}

Moreover, due to
\begin{equation*}
\begin{split}
\bigg(\frac{\widetilde{v}_\xi}{\widetilde{v}}\bigg)_t
&=\bigg(\frac{\widetilde{v}_t}{\widetilde{v}}\bigg)_\xi
=\bigg(\frac{v_t}{v}-\frac{V_t}{V}\bigg)_\xi=\bigg(\frac{s_-v_\xi+u_\xi}{v}-\frac{s_-V_\xi+U_\xi}{V}\bigg)_\xi\\
&=s_-\bigg(\frac{v_\xi}{v}-\frac{V_\xi}{V}\bigg)_\xi-\bigg(\frac{u_\xi}{v}-\frac{U_\xi}{V}\bigg)_\xi
=s_-\bigg(\frac{\widetilde{v}_\xi}{\widetilde{v}}\bigg)_\xi
-\bigg(\frac{\psi_{\xi\xi}}{v}-\frac{U_\xi\phi_\xi}{vV}\bigg)_\xi,
\end{split}
\end{equation*}
we can rewrite $(\ref{2.9})_2$ as
\begin{equation}\label{3.26}
\bigg(\mu\frac{\widetilde{v}_\xi}{\widetilde{v}}-\psi_\xi\bigg)_{t}
-s_{-}\bigg(\mu\frac{\widetilde{v}_\xi}{\widetilde{v}}-\psi_\xi\bigg)_\xi
+(p(v)-p(V))_\xi=0.
\end{equation}
Based on \eqref{3.26}, we can deduce the following result
\begin{Lemma} Under the same assumptions listed in Lemma 3.5, it holds for each $0\leq t \leq T$ that
\begin{equation}\label{3.27}
\begin{split}
&\left\|\left(\frac{\widetilde{v}_\xi}{\widetilde{v}}\right)(t)\right\|^2+
\int^t_0\int^\infty_{0}\frac{\gamma \widetilde{v}^2_\xi}{V^\gamma \widetilde{v}^{\gamma+2}}d\xi d\tau \\
\lesssim & \left\|\left(\frac{\widetilde{v}_{0\xi}}{\widetilde{v_0}},\phi_0,\psi_0,\sqrt{\Phi_0},\psi_{0\xi}\right)
\right\|^2+ \delta^{-1}e^{-c_-\beta}
+|s_{-}| \int^t_0\left(\frac{\widetilde{v}_\xi}{\widetilde{v}}\right)^2(\tau,0)d\tau
\end{split}
\end{equation}
and
\begin{equation}\label{3.28}
 \begin{split}
\int^t_0\int^\infty_{0}\frac{\phi^2_{\xi\xi}}{v}d\xi d\tau
\lesssim &  M^{\gamma+1}\int^t_0\int^\infty_{0}\frac{\widetilde{v}^{2}_{\xi}}{V^\gamma{\widetilde{v}}^{\gamma+2}}d\xi d\tau
+\delta^4m  \int^t_0\int^\infty_{0} \phi^2_\xi d\xi d\tau, \\
\lesssim & M^{\gamma+1}\left(\left\|\left(\frac{\widetilde{v}_{0\xi}}{\widetilde{v_0}},
\phi_0,\psi_0,\sqrt{\Phi_0},\psi_{0\xi}\right)\right\|^2
+\delta^{-1}e^{-c_-\beta}+s_{-} \int^t_0\left(\frac{\widetilde{v}_\xi}{\widetilde{v}}\right)^2(\tau,0)d\tau\right).
\end{split}
\end{equation}

\end{Lemma}
\noindent{\bf Proof}: Multiplying the equation (\ref{3.26}) by $\frac{\widetilde{v}_\xi}{\widetilde{v}}$, we get
\begin{equation}\label{3.29}
\begin{split}
&\left\{\frac{\mu}{2}\left(\frac{\widetilde{v}_\xi}{\widetilde{v}}\right)^2
-\psi_\xi\frac{\widetilde{v}_\xi}{\widetilde{v}}\right\}_t
+\frac{\gamma \widetilde{v}^2_\xi}{V^\gamma \widetilde{v}^{\gamma+2}}
-\left\{s_{-}\frac{\mu}{2}\left(\frac{\widetilde{v}_\xi}{\widetilde{v}}\right)^2
-\psi_\xi \left(s_-\frac{\widetilde{v}_\xi}{\widetilde{v}}+\frac{\psi_{\xi\xi}}{v}-\frac{U_\xi \phi_\xi}{vV}\right)\right\}_\xi \\
=&-\frac{\psi_{\xi\xi}^2}{v}-\frac{U_\xi \phi_\xi \psi_{\xi\xi} }{vV}
+\frac{\gamma V_\xi}{V}(p(V)-p(v))\frac{\widetilde{v}_\xi}{\widetilde{v}}.
\end{split}
\end{equation}

Integrating the above identity with respect to $t$ and $\xi$ over $[0,t]\times{\mathbb{R}}_+$, we can obtain
\begin{equation}\label{3.30}
\begin{split}
&\left\|\left(\frac{\widetilde{v}_\xi}{\widetilde{v}}\right)(t)\right\|^2+
\int^t_0\int^\infty_{0}\frac{\gamma \widetilde{v}^2_\xi}{V^\gamma \widetilde{v}^{\gamma+2}}d\xi d\tau\\
\lesssim & \left\|\frac{\widetilde{v}_{0\xi}}{\widetilde{v_0}}\right\|^2
+\|\psi_{\xi}(t)\|^2+\|\psi_{0\xi}\|^2
+\int^t_0\int^\infty_{0}\frac{\psi^2_{\xi\xi}}{v}d\xi d\tau \\
&+\underbrace{\left|\int^t_0\bigg\{s_{-}\frac{\mu}{2}
\bigg(\frac{\widetilde{v}_\xi}{\widetilde{v}}\bigg)^2
-\psi_\xi \big(s_-\frac{\widetilde{v}_\xi}{\widetilde{v}}+\frac{\psi_{\xi\xi}}{v}-\frac{U_\xi \phi_\xi}{vV}\big)\bigg\}(\tau,0)d\tau\right|}_{I_8} \\
&+\underbrace{\int^t_0\int^\infty_{0}\frac{U^2_\xi \phi^2_\xi }{vV^2}d\xi d\tau}_{I_9}+
\underbrace{\left|\int^t_0\int^\infty_{0}\frac{\gamma V_\xi}{V}(p(V)-p(v))\frac{\widetilde{v}_\xi}{\widetilde{v}}d\xi d\tau\right|}_{I_{10}}.
\end{split}
\end{equation}

As to the estimates on $I_j (j=8,9,10)$, we first bound $I_9$ and $I_{10}$  from \eqref{3.2} as follows
\begin{equation}\label{3.31}
\begin{split}
I_9 \lesssim &  m \delta^4 \int^t_0\int^\infty_{0}\phi^2_\xi d\xi d\tau,\\
I_{10}\lesssim &　\frac{1}{2}\int^t_0\int^\infty_{0}\frac{\gamma \widetilde{v}^2_\xi}{V^\gamma \widetilde{v}^{\gamma+2}}d\xi d\tau +\int^t_0\int^\infty_{0}\frac{\gamma V^2_\xi}{V^2}(p(V)-p(v))^2\widetilde{v}^\gamma  d\xi d\tau \\
\lesssim &　\frac{1}{2}\int^t_0\int^\infty_{0}\frac{\gamma \widetilde{v}^2_\xi}{V^\gamma \widetilde{v}^{\gamma+2}}d\xi d\tau
+\delta^4 \int^t_0\int^\infty_{0} v^{-\gamma-2}\phi^2_\xi  d\xi d\tau \\
\lesssim &　\frac{1}{2}\int^t_0\int^\infty_{0}\frac{\gamma \widetilde{v}^2_\xi}{V^\gamma \widetilde{v}^{\gamma+2}}d\xi d\tau + \delta^4 m^{\gamma+2}\int^t_0\int^\infty_{0}\phi^2_\xi  d\xi d\tau.
\end{split}
\end{equation}

For $I_8$, since $ s_- \phi_\xi=\phi_t-\psi_\xi$, we have
\begin{equation}\label{3.32}
\begin{split}
&s_-\frac{\widetilde{v}_\xi}{\widetilde{v}}\psi_\xi+\frac{\psi_{\xi\xi} \psi_\xi}{v}-\frac{U_\xi \phi_\xi \psi_\xi}{vV} \\
=& s_-\bigg(\frac{\phi_{\xi\xi}}{v}-\frac{V_\xi \phi_\xi}{v V}\bigg)\psi_\xi+\frac{\psi_{\xi\xi} \psi_\xi}{v}-\frac{U_\xi \phi_\xi \psi_\xi}{vV} \\
=& \bigg(\frac{\phi_{t\xi}}{v}-\frac{\psi_{\xi\xi}}{v}\bigg)\psi_\xi-\frac{s_-V_\xi \phi_\xi \psi_\xi}{v V}+\frac{\psi_{\xi\xi} \psi_\xi}{v}-\frac{U_\xi \phi_\xi \psi_\xi}{vV} \\
=& \frac{\phi_{t\xi}\psi_\xi}{v}-\frac{(s_-V_\xi+U_\xi) \phi_\xi \psi_\xi}{v V}
=\frac{\phi_{t\xi}\psi_\xi}{v}-\frac{V_t \phi_\xi \psi_\xi}{v V},
\end{split}
\end{equation}
therefore, from Lemma 3.2, $I_8$ can be controlled by
\begin{equation}\label{3.33}
\begin{split}
{I_8} \lesssim & |s_{-}| \int^t_0\left(\frac{\widetilde{v}_\xi}{\widetilde{v}}\right)^2(\tau,0)d\tau
+\left|\int^t_0\frac{\phi_{t\xi}\psi_\xi}{v}(\tau,0)d\tau\right| +\left|\int^t_0\frac{V_t \phi_\xi \psi_\xi}{vV}(\tau,0)d\tau\right| \\
 \lesssim & |s_{-}| \int^t_0\left(\frac{\widetilde{v}_\xi}{\widetilde{v}}\right)^2(\tau,0)d\tau
 +\delta^2 e^{-c_-\beta}.
\end{split}
\end{equation}

Putting these estimates on $I_j (j=8,9,10)$ into \eqref{3.30},
we have
\begin{equation}\label{3.34}
\begin{split}
&\left\|\left(\frac{\widetilde{v}_\xi}{\widetilde{v}}\right)(t)\right\|^2+
\int^t_0\int^\infty_{0}\frac{\gamma \widetilde{v}^2_\xi}{V^\gamma \widetilde{v}^{\gamma+2}}d\xi d\tau\\
\lesssim & \left\|\frac{\widetilde{v}_{0\xi}}{\widetilde{v_0}}\right\|^2
+\|\psi_{\xi}(t)\|^2+\|\psi_{0\xi}\|^2+\delta^2 e^{-c_-\beta}
+\int^t_0\int^\infty_{0}\frac{\psi^2_{\xi\xi}}{v}d\xi d\tau \\
&+\delta^4 m^{\gamma+2}\int^t_0\int^\infty_{0}\phi^2_\xi  d\xi d\tau
+|s_{-}| \int^t_0\left(\frac{\widetilde{v}_\xi}{\widetilde{v}}\right)^2(\tau,0)d\tau
 \\
\lesssim & \left\|\frac{\widetilde{v}_{0\xi}}{\widetilde{v_0}}\right\|^2+
\left\|\left(\phi_0,\psi_0,\sqrt{\Phi_0},\psi_{0\xi}\right)\right\|^2+  \delta^{-1}e^{-c_-\beta}+|s_{-}| \int^t_0\left(\frac{\widetilde{v}_\xi}{\widetilde{v}}\right)^2(\tau,0)d\tau\\
&+\delta^4 m^{\gamma+3}M^{\gamma+1}\left(\left\|\left(\phi_0,\psi_0,\sqrt{\Phi_0},
\psi_{0\xi}\right)\right\|^2+ \delta^{-1} e^{-c_-\beta}\right).
\end{split}
\end{equation}
Noticing $\delta^2 m^{\gamma+3}M^{\gamma+1}\leq \epsilon_2$,  we can further obtain that
\begin{equation}\label{3.35}
\begin{split}
&\left\|\left(\frac{\widetilde{v}_\xi}{\widetilde{v}}\right)(t)\right\|^2+
\int^t_0\int^\infty_{0}\frac{\gamma \widetilde{v}^2_\xi}{V^\gamma \widetilde{v}^{\gamma+2}}d\xi d\tau\\
\lesssim & \left\|\left(\frac{\widetilde{v}_{0\xi}}{\widetilde{v_0}},
\phi_0,\psi_0,\sqrt{\Phi_0},\psi_{0\xi}\right)\right\|^2+ \delta^{-1}e^{-c_-\beta}
+|s_{-}| \int^t_0\left(\frac{\widetilde{v}_\xi}{\widetilde{v}}\right)^2(\tau,0)d\tau.
\end{split}
\end{equation}
Due to
 \begin{equation}\label{3.36}
\frac{\widetilde{v}_\xi}{\widetilde{v}}=\frac{\phi_{\xi\xi}}{v}-\frac{V_\xi \phi_\xi}{vV},
\end{equation}
thus
\begin{equation}\label{3.37}
\begin{split}
\int^t_0\int^\infty_{0}\frac{\phi^2_{\xi\xi}}{v}d\xi d\tau
\lesssim &  \int^t_0\int^\infty_{0}\left(\frac{v^{\gamma+1}\widetilde{v}^2_{\xi}}{v^\gamma\widetilde{v}^2}
+\frac{V^2_\xi \phi^2_\xi}{vV^2}\right)d\xi d\tau\\
\lesssim &  M^{\gamma+1}\int^t_0\int^\infty_{0}\frac{\widetilde{v}^{2}_{\xi}}{V^\gamma{\widetilde{v}}^{\gamma+2}}d\xi d\tau
+\delta^4m  \int^t_0\int^\infty_{0} \phi^2_\xi d\xi d\tau,
\end{split}
\end{equation}
Combing this result with  \eqref{3.23} and \eqref{3.35},
 one can easily get the estimate (\ref{3.27}). This completes the proof of Lemma 3.6.\\

So far, the only thing left is to control the boundary term  $\int^t_0\big(\frac{\widetilde{v}_\xi}{\widetilde{v}}\big)^2(\tau,0)d\tau$. For this purpose, we can get by recalling the equation (\ref{3.36}) that
\begin{equation}\label{3.38}
\begin{split}
\int^t_0\left(\frac{\widetilde{v}_\xi}{\widetilde{v}}\right)^2(\tau,0)d\tau \lesssim&
\int^t_0\left(\frac{\phi_{\xi\xi}}{v}\right)^2(\tau,0)d\tau
+\int^t_0\left(\frac{V_\xi \phi_\xi}{vV}\right)^2(\tau,0)d\tau \\
\lesssim & \frac{1}{s^2_-v^2_-}\int^t_0\psi^2_{\xi\xi}(\tau,0)d\tau
+\frac{1}{s^2_-v^2_-}\int^t_0\phi^2_{t\xi}(\tau,0)d\tau
+\int^t_0\left(\frac{V_\xi \phi_\xi}{vV}\right)^2(\tau,0)d\tau \\
\lesssim & \frac{1}{s^2_-v^2_-}\int^t_0\psi^2_{\xi\xi}(\tau,0)d\tau +\frac{1}{s^2_-v^2_-}\int^t_0\phi^2_{t\xi}(\tau,0)d\tau
+\int^t_0\left(\frac{V_\xi \phi_\xi}{vV}\right)^2(\tau,0)d\tau \\
\lesssim & \frac{1}{s^2_-v^2_-}\int^t_0\left|\int^\infty_0\frac{d}{d\xi}\psi^2_{\xi\xi}d\xi \right|d\tau +\frac{C\delta^3 e^{-c_-\beta}}{s^2_-}
+\delta^4 e^{-c_-\beta}\\
\lesssim & \frac{M}{s^2_-}\int^t_0\left\|\frac{\psi_{\xi\xi}}{\sqrt{v}}(\tau)\right\|
\left\|\frac{\psi_{\xi\xi\xi}}{\sqrt{v}}(\tau)\right\|d\tau +\frac{\delta^3e^{-c_-\beta}}{s^2_-}.
\end{split}
\end{equation}
Therefore for some $\lambda>0$ to be chosen later, one has
\begin{equation}\label{3.39}
\begin{split}
&\left\|\left(\frac{\widetilde{v}_\xi}{\widetilde{v}}\right)(t)\right\|^2+
\int^t_0\int^\infty_{0}\frac{\gamma \widetilde{v}^2_\xi}{V^\gamma \widetilde{v}^{\gamma+2}}d\xi d\tau \\
\lesssim & \left\|\left(\frac{\widetilde{v}_{0\xi}}{\widetilde{v_0}},\phi_0,
\psi_0,\sqrt{\Phi_0},\psi_{0\xi}\right)\right\|^2+\delta^{-1} e^{-c_-\beta}\\
&+\frac{M}{s_-}\int^t_0\left\|\frac{\psi_{\xi\xi}}{\sqrt{v}}(\tau)\right\|
\left\|\frac{\psi_{\xi\xi\xi}}{\sqrt{v}}(\tau)\right\|d\tau +\frac{C\delta^3 e^{-c_-\beta}}{s_-}\\
\lesssim &  \left\|\left(\frac{\widetilde{v}_{0\xi}}{\widetilde{v_0}},\phi_0,
\psi_0,\sqrt{\Phi_0},\psi_{0\xi}\right)\right\|^2+\delta^{-1} e^{-c_-\beta}\\
&+\frac{M}{\lambda s^2_-}\int^t_0\left\|\frac{\psi_{\xi\xi}}{\sqrt{v}}(\tau)\right\|^2d\tau
+M \lambda\int^t_0\left\|\frac{\psi_{\xi\xi\xi}}{\sqrt{v}}(\tau)\right\|^2d\tau,
\end{split}
\end{equation}
which means that the estimate on the boundary term $\int^t_0\big(\frac{\widetilde{v}_\xi}{\widetilde{v}}\big)^2(\tau,0)d\tau$ can be reduced to the estimate on the term $\int^t_0\|\frac{\psi_{\xi\xi\xi}}{\sqrt{v}}(\tau)\|^2d\tau $, which is the main content of the next lemma.

\begin{Lemma} Under the same assumptions listed in Lemma 3.5, it holds for each $0\leq t\leq T$ that
\begin{equation}\label{3.40}
\begin{split}
&\left\|\psi_{\xi\xi}(t)\right\|^2+\int^t_0\int^\infty_0\frac{\psi^2_{\xi\xi\xi}}{v}d\xi d\tau\\
\lesssim &\left\|\psi_{0\xi\xi}\right\|^2
+m^{2\gamma}M^{\gamma+1}\int^t_0\int^\infty_{0}\frac{\gamma \widetilde{v}^2_\xi}{V^\gamma \widetilde{v}^{\gamma+2}}d\xi d\tau \\
&+ m \sup_{\tau \in[0,t]}\left\{\left\|\frac{\phi_{\xi\xi}}{\sqrt{v}}(\tau)\right\|^2\right\}
\left(\left\|\left(\phi_0,\psi_0,\sqrt{\Phi_0},\psi_{0\xi}\right)\right\|^2+ \delta^{-1}e^{-c_-\beta}\right).
\end{split}
\end{equation}
\end{Lemma}
\noindent{\bf Proof}: Differentiating  $(\ref{3.14})_2$ with respect to $\xi$ twice, it holds
\begin{equation*}
\psi_{\xi\xi t}-s_{-}\psi_{\xi\xi\xi}+(p(v)-p(V))_{\xi\xi}-\mu\left(\frac{\psi_{\xi\xi}}{v}\right)_{\xi\xi}
  =\mu\left(\frac{U_\xi}{v}-\frac{U_\xi}{V}\right)_{\xi\xi},
\end{equation*}
then multiplying the above identity by $\psi_{\xi\xi}$, we have
\begin{equation}\label{3.41}
\begin{split}
&\left(\frac{\psi^2_{\xi\xi}}{2}\right)_t+\mu\frac{\psi^2_{\xi\xi\xi}}{v}
+\left(-s_{-}\frac{\psi^2_{\xi\xi}}{2}+ (p(v)-p(V))_{\xi}\psi_{\xi\xi}
-\mu\left(\frac{\psi_{\xi\xi}}{v}\right)_{\xi}\psi_{\xi\xi}
+\mu\left(\frac{U_\xi\phi_\xi}{vV}\right)_\xi\psi_{\xi\xi}\right)_\xi\\
=&(p(v)-p(V))_{\xi}\psi_{\xi\xi\xi}-\mu\frac{\psi_{\xi\xi}v_\xi}{v^2}\psi_{\xi\xi\xi}
+\mu\left(\frac{U_\xi\phi_\xi}{vV}\right)_{\xi}\psi_{\xi\xi\xi}.
\end{split}
\end{equation}

Integrating the above equation with respect to $t$ and $\xi$ over $[0,t]\times \mathbb{R}_-$, we get
\begin{equation}\label{3.42}
\begin{split}
&\left\|\psi_{\xi\xi}(t)\right\|^2+\int^t_0\int^\infty_0\frac{\psi^2_{\xi\xi\xi}}{v}d\xi d\tau\\
\lesssim &\left\|\psi_{0\xi\xi}\right\|^2
+\underbrace{\left|\int^t_0\left(-s_{-}\frac{\psi^2_{\xi\xi}}{2}+ (p(v)-p(V))_{\xi}\psi_{\xi\xi}
-\mu\left(\frac{\psi_{\xi\xi}}{v}\right)_{\xi}\psi_{\xi\xi}
+\mu\left(\frac{U_\xi\phi_\xi}{vV}\right)_\xi\psi_{\xi\xi}\right)(\tau,0)d\tau\right|}_{I_{11}}\\
&+\underbrace{\left|\int^t_0\int^\infty_0v(p(v)-p(V))^2_{\xi}d\xi d\tau\right|}_{I_{12}}
+\underbrace{\left|\int^t_0\int^\infty_0v\left(\frac{U_\xi\phi_\xi}{vV}\right)^2_{\xi}d\xi d\tau\right|}_{I_{13}}
+\underbrace{\int^t_0\int^\infty_0\frac{\psi^2_{\xi\xi}v^2_\xi}{v^3}d\xi d\tau}_{I_{14}}.
\end{split}
\end{equation}

Now we deal with $I_j (j=11, 12, 13, 14)$ term by term. To do so, since
\begin{equation*}
\begin{split}
&-s_{-}\frac{\psi^2_{\xi\xi}}{2}+ \left\{(p(v)-p(V))_{\xi}
-\mu\left(\frac{\psi_{\xi\xi}}{v}\right)_{\xi}
+\mu\left(\frac{U_\xi\phi_\xi}{vV}\right)_\xi\right\}\psi_{\xi\xi}\\
=& -s_{-}\frac{\psi^2_{\xi\xi}}{2}+\psi_{\xi\xi}(-\psi_{\xi t}+s_{-}\psi_{\xi\xi})
= \frac{s_{-}}{2}\psi^2_{\xi\xi}-\psi_{\xi\xi}\psi_{\xi t}
\end{split}
\end{equation*}
and noticing that
$$
\psi_{t \xi }(t,0)=(s-s_-)U'(-(s-s_-)t+\sigma-\beta),
$$
we can thus bound $I_{11}$ from Lemma 3.2 and Proposition 2.1 as follows
\begin{equation*}
\begin{split}
I_{11}
&\lesssim \left|\frac{s_{-}}{2}\int^t_0\psi^2_{\xi\xi}(\tau,0)d\tau\right|
+\left|\int^t_0\big(\psi_{\xi\xi}\psi_{t\xi}\big)(\tau,0)d\tau\right|\\
&\lesssim \left|\frac{s_{-}}{2}\int^t_0\int^\infty_0\frac{d}{d\xi}\psi^2_{\xi\xi}d\xi
d\tau\right|
+\left|\int^t_0\big(\psi^2_{\xi\xi}\psi_{t\xi}\big)(\tau,0)d\tau\right|
+\left|\int^t_0\psi_{t\xi}(\tau,0)d\tau\right|\\
&\lesssim \left|{s_{-}}\int^t_0\int^\infty_0\psi_{\xi\xi}\psi_{\xi\xi\xi}d\xi
d\tau\right|+\delta^2\left|\int^t_0\psi^2_{\xi\xi}(\tau,0)d\tau\right|
+\delta e^{-c_-\beta} \\
&\lesssim\left(|s_-|+\delta^2\right) \int^t_0\int^\infty_0\frac{\psi^2_{\xi\xi}}{v}d\xi+
\left(|s_-|+\delta^2\right) \int^t_0\int^\infty_0\frac{\psi^2_{\xi\xi\xi}}{v}d\xi
d\tau+ \delta e^{-c_-\beta}.
\end{split}
\end{equation*}
On the other hand, since
\begin{equation}\label{3.43}
\begin{split}
\{(p(v)-p(V))_\xi\}^2&=\{(p'(v)\phi_{\xi\xi}+(p'(v)-p'(V))V_\xi\}^2\\
&=\{(p'(v)\phi_{\xi\xi}+p''(\theta v+(1-\theta)V))V_\xi\phi_\xi\}^2\\
&\leq Cp'(v)^2\phi^2_{\xi\xi}+p''(\theta v+(1-\theta)V))^2V^2_\xi\phi^2_\xi\\
&\lesssim v^{-2\gamma-2}\phi^2_{\xi\xi}+v^{-2\gamma-4}V^2_\xi\phi^2_\xi,
\end{split}
\end{equation}
we can thus control $I_{12}$, $I_{13}$, and $I_{14}$ from \eqref{3.23} as
\begin{equation}\label{3.44}
\begin{split}
I_{12}\lesssim & m^{2\gamma}\int^t_0\int^\infty_0\frac{\phi^2_{\xi\xi}}{v}d\xi d\tau
+m^{2\gamma+3} \delta^4 \int^t_0\int^\infty_0 \phi^2_{\xi}d\xi d\tau \\
\lesssim & m^{2\gamma}\int^t_0\int^\infty_0\frac{\phi^2_{\xi\xi}}{v}d\xi d\tau
+ m^{2\gamma+4}M^{\gamma+1} \delta^4\left(\left\|\left(\phi_0,\psi_0,\sqrt{\Phi_0},\psi_{0\xi}\right)\right\|^2+ \delta^{-1}e^{-c_-\beta}\right),\\
I_{13}\lesssim & \delta^2m\int^t_0\int^\infty_0\frac{\phi^2_{\xi\xi}}{v}d\xi d\tau
+\delta^4 \sup_{\tau \in[0,t]}\Big\{\|\frac{\phi_{\xi}}{v}(\tau)\|^2\Big\}\int^t_0\int^\infty_0\frac{\phi^2_{\xi\xi}}{v}d\xi d\tau
+m^{2} \delta^8 \int^t_0\int^\infty_0 \phi^2_{\xi}d\xi d\tau \\
\lesssim & \delta^4 M^{\gamma-1}\sup_{\tau \in[0,t]}\Big\{\|\Phi(\tau)\|^2\Big\}\int^t_0\int^\infty_0\frac{\phi^2_{\xi\xi}}{v}d\xi d\tau
+m^{2} \delta^8 \int^t_0\int^\infty_0 \phi^2_{\xi}d\xi d\tau \\
\lesssim & \delta^4 M^{\gamma-1}\left(\left\|\left(\phi_0,\psi_0,\sqrt{\Phi_0},\psi_{0\xi}
\right)\right\|^2+\delta^{-1}e^{-c_-\beta} \right)\int^t_0\int^\infty_0\frac{\phi^2_{\xi\xi}}{v}d\xi d\tau,
\end{split}
\end{equation}
and
\begin{equation}\label{3.45}
\begin{split}
I_{14}\lesssim &  \int^t_0\int^\infty_0\frac{\psi^2_{\xi\xi}\left(\phi^2_{\xi\xi}+V^2_{\xi}\right)}{v^3}d\xi d\tau \\
\lesssim &  m \sup_{\tau \in[0,t]}\left\{\left\|\frac{\phi_{\xi\xi}}{\sqrt{v}}(\tau)\right\|^2\right\}\int^t_0\int^\infty_0\frac{\psi^2_{\xi\xi}}{v}d\xi d\tau \\
\lesssim & m \sup_{\tau \in[0,t]}\left\{\left\|\frac{\phi_{\xi\xi}}{\sqrt{v}}(\tau)\right\|^2\right\}
\left(\left\|\left(\phi_0,\psi_0,\sqrt{\Phi_0},\psi_{0\xi}\right)\right\|^2+ \delta^{-1}e^{-c_-\beta}\right).
\end{split}
\end{equation}

Inserting these estimates into \eqref{3.42},  we finally arrive at
\begin{equation}\label{3.46}
\begin{split}
&\|\psi_{\xi\xi}(t)\|^2+\int^t_0\int^\infty_0\frac{\psi^2_{\xi\xi\xi}}{v}d\xi d\tau\\
\lesssim &\|\psi_{0\xi\xi}\|^2+m^{2\gamma}\int^t_0\int^\infty_0\frac{\phi^2_{\xi\xi}}{v}d\xi d\tau
+ m^{2\gamma+4}M^{\gamma+1} \delta^4\left(\left\|\left(\phi_0,\psi_0,\sqrt{\Phi_0},\psi_{0\xi}\right)\right\|^2+ \delta^{-1}e^{-c_-\beta}\right)\\
&+  m \sup_{\tau \in[0,t]}\left\{\left\|\frac{\phi_{\xi\xi}}{\sqrt{v}}(\tau)\right\|^2\right\}
\left(\left\|\left(\phi_0,\psi_0,\sqrt{\Phi_0},\psi_{0\xi}\right)\right\|^2+ \delta^{-1}e^{-c_-\beta}\right).
\end{split}
\end{equation}

Recalling that $ m^{2\gamma+4}M^{\gamma+1} \delta^4\ll 1$ and \eqref{3.28}, one can thus get the estimate \eqref{3.40} immediately and the proof of Lemma 3.7 is complete.\\

Setting
\begin{equation}\label{3.47}
H(t):=\left\|\frac{\widetilde{v}_\xi}{\widetilde{v}}(t)\right\|^2+
\int^t_0\int^\infty_{0}\frac{\gamma \widetilde{v}^2_\xi}{V^\gamma \widetilde{v}^{\gamma+2}}d\xi d\tau,
\end{equation}
we get by combing (\ref{3.36}) and (\ref{3.24}) together that
\begin{equation}\label{3.48}
\begin{split}
\left\|\frac{\phi_{\xi\xi}}{\sqrt{v}}(t)\right\|^2
&\lesssim M\left\|\frac{\widetilde{v}_\xi}{\widetilde{v}}(t)
\right\|^2+m \delta^4 \|\phi_{\xi}(t)\|^2\\
&\lesssim M H(t)+m M^{\gamma+1}\delta^4\left(\left\|\left(\phi_0,\psi_0,\sqrt{\Phi_0},\psi_{0\xi}\right)\right\|^2+ \delta^{-1}e^{-c_-\beta}\right),
\end{split}
\end{equation}
from which and the estimate (\ref{3.40}), one can get that
\begin{equation}\label{3.49}
\int^t_0\int^\infty_0\frac{\psi^2_{\xi\xi\xi}}{v}d\xi d\tau
\lesssim\|\psi_{0\xi\xi}\|^2+m^{2\gamma}M^{\gamma+1}H(t).
\end{equation}

Combining \eqref{3.24}, \eqref{3.48}, and \eqref{3.49} with  (\ref{3.39}), we see that
\begin{equation}\label{3.50}
\begin{split}
H(t)\lesssim& \left\|\left(\frac{\widetilde{v}_{0\xi}}{\widetilde{v_0}},\phi_0,\psi_0,
\sqrt{\Phi_0},\psi_{0\xi}\right)\right\|^2+\delta^{-1} e^{-c_-\beta}\\
&+\frac{ M}{\lambda s^2_-}\left(\left\|\left(\phi_0,\psi_0,\sqrt{\Phi_0},\psi_{0\xi}\right)\right\|^2
+\delta^{-1} e^{-c_-\beta}\right)\\
&+\lambda M \left(\left\|\psi_{0\xi\xi}\right\|^2+m^{2\gamma}M^{\gamma+1}H(t)\right)\\
\lesssim& \left\|\left(\frac{\widetilde{v}_{0\xi}}{\widetilde{v_0}},\phi_0,
\psi_0,\sqrt{\Phi_0},\psi_{0\xi}\right)\right\|^2
+\frac{M}{\lambda s^2_-}\left(\left\|\left(\phi_0,\psi_0,\sqrt{\Phi_0},\psi_{0\xi}\right)\right\|^2
+\delta^{-1} e^{-c_-\beta}\right)\\
& +\lambda M\left\|\psi_{0\xi\xi}\right\|^2+\lambda m^{2\gamma}M^{\gamma+2}H(t).
\end{split}
\end{equation}

Having obtained \eqref{3.50}, if we let $\lambda\sim N(t)^{\frac{1}{2}}$ and assume further that
\begin{equation}\label{3.51}
N(t)^{\frac{1}{2}}  m^{2\gamma}M^{\gamma+2}< 1,
\end{equation}
we can obtain that
\begin{equation}\label{3.52}
H(t)\lesssim \left\|\left(\frac{\widetilde{v}_{0\xi}}{\widetilde{v_0}},\psi_{0\xi\xi}
\right)\right\|^2
+\frac{M}{\lambda s^2_-}\left(\left\|\left(\phi_0,\psi_0,\sqrt{\Phi_0},\psi_{0\xi}\right)\right\|^2
+\delta^{-1} e^{-c_-\beta}\right).
\end{equation}
From the above estimate, we can get that
\begin{Lemma} Under the same assumptions listed in Lemma 3.5, if we assume further that
\begin{equation}\label{3.53}
\begin{cases}
N(t) m^{4\gamma}M^{2\gamma+4}< 1,\\
M\left(\left\|\left(\phi_0,\psi_0,\sqrt{\Phi_0},\psi_{0\xi}\right)\right\|^2
+\delta^{-1}e^{-c_-\beta}\right)\leq
N(t)^\frac{1}{2} s^2_- \left\|\left(\frac{\widetilde{v}_{0\xi}}{\widetilde{v_0}},\psi_{0\xi\xi}\right)
\right\|^2,
\end{cases}
\end{equation}
then it holds for each $0\leq t \leq T$ that
\begin{equation}\label{3.54}
H(t)\lesssim  \left\|\left(\frac{\widetilde{v}_{0\xi}}{\widetilde{v_0}},
\psi_{0\xi\xi},\phi_0,\psi_0,\sqrt{\Phi_0},\psi_{0\xi}\right)\right\|^2
+\delta^{-1}e^{-c_-\beta}.
\end{equation}
\end{Lemma}

In summary, what we have obtained up to now is that if
\begin{equation}\label{3.55}
\begin{cases}
N(t)m^{2\gamma}M^{2\gamma+4}\ll 1,,\\
\delta^2 m^{\gamma+3}M^{\gamma+1}\ll 1,\\
s^2_- mM^{\gamma+1} \ll 1, \\
M\left(\left\|\left(\phi_0,\psi_0,\sqrt{\Phi_0},\psi_{0\xi}
\right)\right\|^2+\delta^{-1}e^{-c_-\beta}\right)\leq
N(t)^\frac{1}{2} s^2_-\left\|\left(\frac{\widetilde{v}_{0\xi}}{\widetilde{v_0}},
\psi_{0\xi\xi}\right)\right\|^2,
\end{cases}
\end{equation}
then it holds for each $0\leq t \leq T$ that
\begin{equation}\label{3.56}
\begin{split}
&\left\|\left(\phi,\psi,\sqrt{\Phi},\psi_\xi\right)(t)\right\|^2
+ \int^t_0\int^\infty_{0}\left({\psi}^2_{\xi}+ \frac{{\psi}^2_{\xi\xi}}{v}\right)d\xi d\tau
\lesssim \left\|\left(\phi_0,\psi_0,\sqrt{\Phi_0},\psi_{0\xi}\right)\right\|^2
+\delta^{-1}e^{-c_-\beta},\\
&\left\|\frac{\widetilde{v}_\xi}{\widetilde{v}}(t)\right\|^2+
\int^t_0\int^\infty_{0}\frac{\gamma \widetilde{v}^2_\xi}{V^\gamma \widetilde{v}^{\gamma+2}}d\xi d\tau
\lesssim \left\|\left(\frac{\widetilde{v}_{0\xi}}{\widetilde{v_0}},\psi_{0\xi\xi},
\phi_0,\psi_0,\sqrt{\Phi_0},\psi_{0\xi}\right)\right\|^2
+\delta^{-1}e^{-c_-\beta}.
\end{split}
\end{equation}

Having obtained \eqref{3.56}, we now turn to deduce the desired uniform positive lower and upper bounds on $v(t,x)$ in terms of the initial perturbation. In fact from \eqref{3.56}, we can get by employing Y. Kanel's argument \cite{Kanel-Diff Eqns-1968} as in \cite{Fan-Liu-Wang-Zhao-JDE-2014} that
\begin{Lemma} Under the assumption (\ref{3.55}), we can get that
\begin{equation}\label{3.57}
C^{-1}B^\frac{2}{1-\gamma}_0\leq v(t,\xi) \leq C B^2_0,\quad \forall (t,\xi)\in [0,T]\times\mathbb{R}_+.
\end{equation}
Here
\begin{equation*}
 B_0:=\left(\left\|\left(\phi_0,\psi_0,\sqrt{\Phi_0},\psi_{0\xi}\right)\right\|^2
 +\delta^{-1}e^{-c_-\beta}\right)^\frac{1}{2}
\left(\left\|\left(\frac{\widetilde{v}_{0\xi}}{\widetilde{v_0}},
\psi_{0\xi\xi},\phi_0,\psi_0,\sqrt{\Phi_0},\psi_{0\xi}\right)\right\|^2
+\delta^{-1}e^{-c_-\beta}\right)^\frac{1}{2}.
\end{equation*}
Moreover, as a direct consequence of \eqref{3.56} and \eqref{3.57}, we can deduce that
\begin{equation*}
\begin{split}
&\left\|\phi_{\xi\xi}(t)\right\|^2+\int^t_0\left\|\phi_{\xi\xi}(\tau)\right\|^2 d\tau
\leq C\left(\left\|\left(\phi_0,\psi_0\right)\right\|_2,\delta^{-1}e^{-c_-\beta}\right),\\
&\left\|\psi_{\xi\xi}(t)\right\|^2+\int^t_0\left\|\frac{\psi_{\xi\xi\xi}}{\sqrt{v}}(\tau)
\right\|^2 d\tau
\leq C\left(\left\|\left(\phi_0,\psi_0\right)\right\|_2,\delta^{-1}e^{-c_-\beta}\right).
\end{split}
\end{equation*}
Here $C\left(\left\|\left(\phi_0,\psi_0\right)\right\|_2,\delta^{-1}e^{-c_-\beta}\right)$
is some positive constant depending only on $\left\|\left(\phi_0,\psi_0\right)\right\|_2$ and $\delta^{-1}e^{-c_-\beta}$.
\end{Lemma}

\subsection{The proof of our main result}
With the above preparations in hand, we now turn to prove Theorem 1. Noticing that
$$
\Phi_0(x)\lesssim \left(|V(0,x)|^{-\gamma-1}+|v(0,x)|^{-\gamma-1}\right)\phi^2_{0x},
$$
we first deduce from the assumptions (H$_1$)-(H$_3$) that
\begin{equation}\label{3.58}
\begin{split}
&\left\|\sqrt{\Phi_0}\right\|\lesssim \left(1+\delta^{-\frac{(\gamma+1)l}{2}}\right)\left\|\phi_{0x}\right\|, \\
&\left\|\frac{\widetilde{v}_{0\xi}}{\widetilde{v_0}}\right\|\lesssim \delta^{-l}\left(\left\|\phi_{0xx}\right\|+\delta^2\left\|\phi_{0x}\right\|\right).
\end{split}
\end{equation}
Hence, if we ask $\beta=o(\delta^{-1})$ and noticing that $c_-=O(1)\delta$, one can deduce that
\begin{equation}\label{3.59}
\delta^{-1}e^{-c_-\beta}\leq\delta^{2\alpha-(\gamma+1)l}
\end{equation}
provided that $\delta>0$ satisfies $0<\delta\leq \delta_1$ holds for some sufficiently small $\delta_1>0$.

Consequently
\begin{equation}\label{3.60}
\begin{split}
&\left\|\left(\phi_0,\psi_0,\sqrt{\Phi_0},\psi_{0\xi}\right)\right\|
+\delta^{-1}e^{-c_- \beta}
\lesssim \delta^{\alpha-\frac{(\gamma+1)l}{2}}+\delta^{-1}e^{-c_- \beta}
\lesssim \delta^{\alpha-\frac{(\gamma+1)l}{2}},\\
&\left\|\left(\psi_{0\xi\xi},\frac{\widetilde{v}_{0\xi}}{\widetilde{v_0}}
\right)\right\| \lesssim\delta^{-\kappa-l}.
\end{split}
\end{equation}

Now we prove Theorem 1 by exploiting the continuation argument. Applying Proposition \ref{locals}, we can find a positive
constant $t_0$, which only depends on $\delta$ and $\|(\phi_0,\psi_0)\|_2$ such that the  problem \eqref{2.9} admits
a unique solution $(\phi(t,\xi),\psi(t,\xi))\in X_{m_0,M_0}(0,t_0)$ with
\begin{equation*}
m_0=2^{-1}C^{-1}_0\delta^l, \quad M_0=2C_0\left(1+\delta^{-l}\right),
\end{equation*}
and we have from Sobolev's inequality that for each $0\leq t\leq t_0$
\begin{equation*}
\sup_{[0,t_0]}\Big\{\|(\phi,\psi)(t)\|_{L^\infty}\Big\}\leq \|(\phi,\psi)(t)\|_1
\leq C\|(\phi_0,\psi_0)\|_1 \leq C\delta^\alpha.
\end{equation*}
Thus we can take $N(t_0)=C\delta^\alpha$.

If
\begin{equation*}
\begin{cases}
&N(t_0)m_0^{-4\gamma}M_0^{2\gamma+4}\ll 1,\\
&\delta^2 m^{-\gamma-3}_0M^{\gamma+1}_0\ll 1,\\
&s^2_- m_0^{-1}M^{-\gamma-1}_0 \ll 1, \\
&M_0\left(\left\|\left(\phi_0,\psi_0,\sqrt{\Phi_0},\psi_{0\xi}\right)\right\|^2
+\delta^{-1}e^{-c_-\beta}\right)\leq
N(t_0)^\frac{1}{2} s^2_-\left\|\left(\frac{\widetilde{v}_{0\xi}}{\widetilde{v_0}},
\psi_{0\xi\xi}\right)\right\|^2,
\end{cases}
\end{equation*}
which is equivalent to ask that
\begin{equation}\label{alpha-beta-1}
\begin{cases}
\delta^{\alpha-(6\gamma+4)l}<1,\\
\delta^{2-2(\gamma+2)l}<1,\\
\delta^{2h-(\gamma+2)l}<1,\\
\delta^{2\alpha-(\gamma+1)l-l-\frac{\alpha}{2}-2h}\leq \delta^{-2\kappa-2l},
\end{cases}
\end{equation}
then we can deduce that the estimates obtained in Lemma 3.2-lemma 3.8 hold with $0\leq t \leq t_0$, $m=m^{-1}_0$ and $M=M_0$. Thus from the result of Lemma 3.10, it yields that
\begin{equation}\label{3.62}
2^{-1}C_3^{-1}\delta^{ \frac{2}{1-\gamma} \left\{\left(\alpha-\frac{(\gamma+1)l}{2}\right)-(\kappa+l)\right\} }
\leq  v(t,\xi) \leq 2C_3\delta^{ 2\left\{\left(\alpha-\frac{(\gamma+1)l}{2}\right)-(\kappa+l)\right\} }
\end{equation}
holds for each $0\leq t \leq t_0$ and
\begin{equation}\label{3.63}
\left\|\left(\phi,\psi,\sqrt{\Phi},\psi_{\xi}\right)(t)\right\|\lesssim \delta^{\alpha-\frac{(\gamma+1)l}{2}}
\end{equation}
and
\begin{equation}\label{3.64}
\|(\phi,\psi)(t)\|^2_2+\int^{t}_0\left(\|\phi_\xi(\tau)\|^2_1
+\|\psi_{\xi}(\tau)\|^2_2\right)d\tau
\lesssim C\left(\left\|(\phi_0,\psi_0)\right\|_2,\delta\right)
\end{equation}
hold for each $0\leq t \leq t_0$.

Note that if the parameters $l, h, \alpha,$ and $\kappa$ satisfy
\begin{equation}\label{3.65}
\begin{cases}
1>(\gamma+2)l,\quad \alpha> (6\gamma+4)l,\\[2mm]
\frac{(\gamma+2)l}{2}<h<\kappa+\frac{3\alpha}{4}-\frac{\gamma}{2}l,
\end{cases}
\end{equation}
then there exists a sufficiently small $0<\delta_2\leq\delta_1$ such that for all $\delta\in (0,\delta_2]$, \eqref{alpha-beta-1} hold.

Next if we take $(\phi(t_1,\xi),\psi(t_1,\xi))$ as the initial data,
we can deduce by employing Proposition \ref{locals} again that the unique local solution
$(\phi(t,\xi),\psi(t,\xi))$ constructed above can be extended to the time internal $[t_0,t_0+t_1]$ and satisfies
\begin{equation*}
\sup_{[0,t_0+t_1]}\Big\{\|(\phi,\psi)(t)\|_{L^\infty}\Big\}\leq \max\{N(t_0), 2\|(\phi,\psi)(t_0)\|_1\}\leq C_4\delta^{\alpha-\frac{(\gamma+1)l}{2}},
\end{equation*}
and
$$
\frac{m_1}{2}\leq v(t,\xi)\leq 2M_1,\quad \forall (t,\xi)\in[0,t_1+t_2]\times\mathbb{R}_+
$$
with
\begin{equation*}
\begin{split}
& m_1:=2^{-1}C_3^{-1} \delta^{ \frac{2}{1-\gamma} \left\{\alpha-\frac{(\gamma+1)l}{2}-(\kappa+l)\right\}},\\
& M_1:=2C_3\delta^{ 2\left\{\alpha-\frac{(\gamma+1)l}{2}-(\kappa+l)\right\} }.
\end{split}
\end{equation*}
Thus we can take $N(t_0+t_1)=C_4\delta^{\alpha-\frac{(\gamma+1)l}{2}}$ and if we assume that
\begin{equation*}
\begin{cases}
&N(t_0+t_1)m_1^{-4\gamma}M_1^{2\gamma+4}\ll 1,\\
&\delta^2 m^{-\gamma-3}_1M^{\gamma+1}_1\ll 1,\\
&s^2_- m_1^{-1}M^{-\gamma-1}_1 \ll 1, \\
&M_1\left(\left\|\left(\phi_0,\psi_0,\sqrt{\Phi_0},\psi_{0\xi}\right)\right\|^2
+\delta^{-1}e^{-c_-\beta}\right)\leq
N(t_0+t_1)^\frac{1}{2} s^2_- \left\|\left(\frac{\widetilde{v}_{0\xi}}{\widetilde{v_0}},\psi_{0\xi\xi}\right)\right\|^2,
\end{cases}
\end{equation*}
that is
\begin{equation*}
\begin{cases}
\delta^{\alpha-\frac{(\gamma+1)l}{2}
-\frac{4(\gamma^2+3\gamma-2)}{\gamma-1}\left(\kappa+l-
\left(\alpha-\frac{(\gamma+1)l}{2}\right)\right)}<1,\\
\delta^{2
-\frac{2\gamma^2+2\gamma+4}{\gamma-1}\left(\kappa+l-
\left(\alpha-\frac{(\gamma+1)l}{2}\right)\right)}<1,\\
\delta^{2h-\frac{2\gamma^2}{\gamma-1}\left(\kappa+l-
\left(\alpha-\frac{(\gamma+1)l}{2}\right)\right)}<1,\\
\delta^{{2\left(\alpha-\frac{(\gamma+1)l}{2}\right)
-2\left(\kappa+l-
\left(\alpha-\frac{(\gamma+1)l}{2}\right)\right)}}
\leq \delta^{\frac{1}{2}\left(\alpha-\frac{(\gamma+1)l}{2}\right)
+2h-2\kappa-2l},
\end{cases}
\end{equation*}
which hold for all $\delta\in (0,\delta_3]$ for some sufficiently small constant $0<\delta_3\leq \delta_1$ provided that
\begin{equation}\label{3.66}
\begin{cases}
0<\theta < \frac{\gamma-1}{4(\gamma^2+3\gamma-2)}\left(\alpha-\frac{(\gamma+1)l}{2}
\right),\\
 0<\theta < \frac{\gamma-1}{\gamma^2+\gamma+2},\\
0<\theta<\frac{\gamma-1}{\gamma^2}h,\\
0<h <\frac{7}{4}\left(\alpha-\frac{(\gamma+1)l}{2}\right),
\end{cases}
\end{equation}
where
$$
\theta:=\kappa+l-\left(\alpha-\frac{(\gamma+1)l}{2}\right),
$$
then the assumptions listed in Lemma 3.3-Lemma 3.8 hold with $t_0\leq t \leq t_0+t_1$, $m=m^{-1}_1$ and $M=M_1$ and consequently the estimates (\ref{3.55}), (\ref{3.56}) and (\ref{3.57}) obtained in Lemma 3.2-Lemma 3.8 hold for each  $0 \leq t \leq t_0+t_1$. If we take
$(\phi( t_0+t_1,\xi),\psi( t_0+t_1,\xi))$ as the initial data and employ Proposition \ref{locals} again, we can then
extend the above solution $(\phi( t_0+t_1,\xi),\psi(t,\xi))$ to the time step $t=t_0+2t_1$.
Repeating the above procedure, we thus extend $(\phi(t,\xi),\psi(t,\xi))$ step by step to the unique
global solution and the estimates (\ref{3.55}), (\ref{3.56}) and (\ref{3.57}) hold for each  $ t \geq 0$ provided that the conditions (H$_1$)-(H$_3$) are assumed to be hold, the initial perturbation $(\phi_0(\xi),\psi_0(\xi))$ satisfies \eqref{2.13} with the parameters $\alpha, \kappa$ and $l$ satisfy \eqref{3.65}-\eqref{3.66}, and $\delta$ satisfies $0<\delta\leq\delta_0:=\min\{\delta_2,\delta_3\}$. Noticing that \eqref{3.65} and \eqref{3.66} are nothing but the assumption \eqref{2.14} imposed in Theorem 1, we thus  complete the proof of Theorem 1.

\bigbreak

\section{Acknowledgement}
Dongfen Bian is partially supported by two grants of the National Natural Science Foundation of China under the contracts 11471323 and 11271052 respectively and by a grant of the China Postdoctoral Science Foundation under contract 2015M570939. Lili Fan is partially supported by a grant from the National Natural Science Foundation of China under contract 11301405. Lin He is partially supported by  ``the Fundamental Research Funds for the Central Universities". Huijiang Zhao is partially supported by three grants from the National Natural Science Foundation of China under contracts 10925103, 11271160, and 11261160485, respectively.

\end{document}